\documentclass[pdflatex,sn-mathphys]{sn-jnl}
\jyear{2021}%

\theoremstyle{thmstyleone}%
\theoremstyle{thmstyletwo}%

\theoremstyle{thmstylethree}%

\usepackage{setspace}
\usepackage{color,soul}
\usepackage{textcomp}
\usepackage{gensymb}
\usepackage{amsmath,amssymb,bm}

\usepackage{relsize}
\usepackage{nomencl}
\usepackage{soul}
\usepackage{graphicx,subfig,float}
\usepackage{geometry}
\usepackage{pdflscape,rotating}
\usepackage[normalem]{ulem}
\usepackage{lineno,hyperref,cleveref}
\usepackage{lipsum}
\usepackage{multirow}
\usepackage{mathtools}

\makenomenclature

\begin{document}

\title[Data Assimilation with Noisy Observations for the Rayleigh-B\'enard Convection]{Continuous and Discrete Data Assimilation with Noisy Observations for the Rayleigh-B\'enard Convection: A Computational Study}

\author[1]{\fnm{Mohamad Abed El Rahman} \sur{Hammoud}}
\author[2]{\fnm{Olivier} \sur{Le Ma\^{\i}tre}}
\author[3,4]{\fnm{Edriss S.} \sur{Titi}}
\author[1]{\fnm{Ibrahim} \sur{Hoteit}}
\author*[5]{\fnm{Omar} \sur{Knio}}
\email{omar.knio@kaust.edu.sa}

\affil[1]{\orgdiv{Division of Physical Sciences and Engineering}, \orgname{King Abdullah University of Science and Technology}, \orgaddress{ \city{Thuwal}, \postcode{23955}, \country{KSA}}}
\affil[2]{\orgname{Centre de Mathématiques Appliquées, CNRS and Inria, Ecole Polytechnique}, \orgaddress{\postcode{91120}, \state{Palaiseau}, \country{France}}}
\affil[3]{\orgdiv{Department of Applied Mathematics and Theoretical Physics}, \orgname{University of Cambridge}, \orgaddress{\city{Cambridge CB3 0WA}, \country{UK}}}
\affil[4]{\orgdiv{Department of Mathematics}, \orgname{Texas A \& M University}, \orgaddress{ \city{College Station}, \postcode{77843}, \state{Texas}, \country{USA}}}
\affil*[5]{\orgdiv{Division of Computer, Electrical and Mathematical Sciences and Engineering}, \orgname{King Abdullah University of Science and Technology}, \orgaddress{ \city{Thuwal}, \postcode{23955}, \country{KSA}}}

\abstract{Obtaining accurate high-resolution representations of model outputs is essential to describe the system dynamics.
In general, however, only spatially- and temporally-coarse observations of the system states are available.
These observations can also be corrupted by noise.
Downscaling  is a process/scheme in which one uses coarse scale observations to reconstruct the high-resolution solution of the system states. Continuous Data Assimilation (CDA) is a recently introduced downscaling algorithm that constructs an increasingly accurate representation of the system states by continuously nudging the large scales using the coarse observations. We introduce a Discrete Data Assimilation (DDA) algorithm as a downscaling algorithm based on CDA with discrete-in-time nudging.
We then investigate the performance of the CDA and DDA algorithms for downscaling noisy observations of the Rayleigh-B\'enard convection system in the chaotic regime. 
In this computational study, a set of noisy observations was generated by perturbing a reference solution with Gaussian noise before downscaling them.
The downscaled fields are then assessed using various error- and ensemble-based skill scores. The CDA solution was shown to converge towards the reference solution faster than that of DDA but at the cost of a higher asymptotic error. 
The numerical results also suggest a quadratic relationship between the $\ell_2$ error and the noise level for both CDA and DDA. 
Cubic and quadratic dependences of the DDA and CDA expected errors on the spatial resolution of the observations were obtained, respectively.}


\maketitle

\section{Introduction}
\label{sec:intro}

High-resolution representations of the states of dynamical systems are essential in various fields \citep{Moser1987}. 
Downscaling techniques were introduced to synchronize the fine-scale features of a dynamical system to a reference solution represented by coarse-scale observations \citep{Wilby1997}. 
These techniques were applied in various applications, such as oceanic \citep{KATAVOUTA2016, Sang2020} and atmospheric modeling \citep{Dosio2015, Srinivas2019}. 
Downscaling techniques can be categorized into two main approaches: statistical downscaling and dynamical downscaling.

Statistical downscaling establishes a statistical relationship between the observed variables at the coarse-scales and the variables of interest at the fine-scales, based on which, the fine-scales can be estimated \citep{Huth2002, TISSEUIL2010, Jha2015, LAFLAMME2016}. 
Dynamical downcaling methods, on the other hand, constrain the fine-scale model solution to the coarse-scale observations. 
Nudging techniques are commonly used for dynamical downscaling, where the fine-scale solution is forced towards the coarse-scale data, point-by-point as in grid nudging and low-frequency to low-frequency as in spectral nudging \citep{vonStorch2000}. 
Several studies compared the performances of statistical and dynamical downscaling methods and reported that both techniques generally perform comparably \citep{Kidson1998, Altaf2017, LEROUX2018, Srinivas2019}.

Continuous Data Assimilation (CDA) is a dynamical downscaling algorithm that constrains the large-scale solution of a system state to coarse-scale observations \citep{Azouani2013}. 
The CDA algorithm introduces a nudging term, as a source term in the equations of motion, to compare the coarse scales of the system state to those of the observations.
The theoretical grounds of CDA state that the downscaled solution has an exponentially increasing accuracy with time, regardless of the initial conditions \citep{Azouani2013}.
The downscaled solution was also proven to asymptote to the true solution, provided sufficient observations in time.
The CDA algorithm was successfully applied to the $2$D Rayleigh-B\'enard equations \citep{FARHAT2015, Altaf2017, Farhat2018}, $3$D Rayleigh-B\'enard equations \citep{FARHAT2016_3dBenard}, quasi-geostrophic equations \citep{Jolly2019}, $2$D Navier Stokes equations \citep{gesho2016}, $3$D Primitive equations \citep{Pei2019} and Global Circulation Models (GCMs) \citep{Srinivas2019}. 
The theoretical framework of CDA with noisy stochastic measurements was more recently established by \cite{Bessaih2015}, who showed that the algorithm asymptotes to the true solution with a tolerance limit proportional to the trace of the covariance of the noisy observations.

This work introduces Discrete Data Assimilation (DDA), a downscaling algorithm that is a discrete-in-time counterpart of CDA, arising in spirit of the works of Hayden et al.~\cite{HAYDEN2011} and Celik et al.~\cite{Celik2019},  where downscaling is performed based on an ``implusive'' update of the state field, using Fourier-based spectral projections specified in terms of the coarse-scale observations.   In contrast, in the present work DDA introduces a source term into the governing equations that is concentrated at the observation times.  Thus, whereas CDA continuously nudges the governing equations with the coarse-scale observations of the system state, DDA only nudges the system at the discrete times when observational data are available.  To numerically analyze the performances of both downscaling algorithms for the case of noisy coarse-scale observations, we rely on synthetic observational data generated by perturbing a reference solution of the Rayleigh-B\'enard convection system in the turbulent regime.
These noisy coarse-scale observations of the velocity and temperature are downscaled realization-by-realization using CDA and DDA, which results in a set of downscaled solutions.

The performances of both CDA and DDA were then assessed by analyzing various skill scores.
Numerical experiments were also conducted to examine the sensitivity of DDA and CDA to the noise level and temporal and spatial resolutions of the observations. 
The statistics of the downscaled ensemble of solutions were then examined and the precision of the skill scores estimators were analyzed 
using a bootstrap analysis. 
Finally, because velocity observations alone are sufficient to fully describe the downscaled Rayleigh-B\'enard solution fields \citep{FARHAT2015}, the critical noise level at which noisy temperature observations should be discarded in favor of employing noisy velocity observations alone was also investigated. 
The numerical results suggest that the $\ell_2$-errors resulting from CDA and DDA are proportional to the observational noise variance and have a cubic relation to the spatial resolution. 
The CDA downscaled fields are also shown to converge to the asymptotic behavior faster than DDA, although it converges
to a lower asymptotic error level in comparison to CDA.

The following section provides an overview of the Rayleigh-B\'enard convection system, the CDA and DDA algorithms, and the skill scores to assess their performances.
Section \ref{sec:expSet} describes the experimental setting, including the generation of noisy observational data, and the model and downscaling parameters. 
The results of the numerical experiments are then presented in Section \ref{sec:results}.
The main conclusions of the study are outlined in Section \ref{sec:conclusion}.

\section{Preliminaries}

\subsection{Rayleigh-B\'enard Convection System}
\label{sec:rb2d}

The Rayleigh-B\'enard convection is a thermo-fluidic instability driven by the temperature difference between a hotter bottom boundary and a cooler top boundary. 
At high Ra, chaos dominates the flow with no particular structures to be observed \citep{Curry1984, Paul2007}. 
The system under study consists of a $2$D rectangular periodic channel of vertical height $L_y$ and horizontal period $L_x$.
The $2$D Boussinesq equations are then recast into their non-dimensional form, using the same scaling used in \cite{LEQUERE1991,LEMAITRE2002}, such that the equations are solved in a fundamental spatial domain $\Omega = \left[0, L_x/L_y \right] \times \left[0, 1 \right]$.
The solution of the governing equations are the dimensionless velocity vector $\mathbf{u} = \left( u, v \right)$, pressure $p$ and temperature $\Theta$.
These equations are expressed as:
\begin{equation}
    \nabla \cdot \mathbf{u} = 0,
    \label{eqn:incompr_Continuity}
\end{equation}
\begin{equation}
    \frac{\partial \mathbf{u}}{\partial t} +  \left( \mathbf{u} \cdot \nabla \right) \mathbf{u} + \nabla p = \frac{Pr}{\sqrt{Ra}} \nabla ^2 \mathbf{u} + Pr \Theta \mathbf{e}_y,
    \label{eqn:momentum_true}
\end{equation}
\begin{equation}
    \frac{\partial \Theta}{\partial t} + \left( \mathbf{u} \cdot \nabla \right) \Theta = \frac{1}{\sqrt{Ra}} \nabla ^2 \Theta + \mathbf{u} \cdot \mathbf{e}_y,
    \label{eqn:energy_true}
\end{equation}
where $Ra \coloneqq \left(g \alpha \Delta T L_y^3\right) \left(\nu \kappa\right)^{-1}$ is the Rayleigh number, $g$ is the gravitational acceleration, $\alpha$ is the thermal expansion coefficient of the fluid, $\Delta T$ is the dimensional temperature difference between the top and bottom boundaries, $\nu$ is the kinematic viscosity, $\kappa$ is the thermal diffusivity coefficient, $Pr = \nu \kappa ^{-1}$ is the Prandtl number, $t$ is time and $\mathbf{e}_y$ is the unit basis vector in vertical direction.

This system of equations can be solved once the initial and boundary conditions are specified.
The temperature and velocities are initialized with independent, spatially uncorrelated random noise sampled point-by-point from uniform distributions $ \mathcal{U}\left(-0.2, 0.2\right) $ and $\mathcal{U}\left(-0.1, 0.1\right)$, respectively, where $\mathcal{U}\left(a,b\right)$ refers to a uniform distribution on the interval $(a,b)$.
Periodic boundary conditions were specified in $x$, whereas at the horizontal boundaries, no-slip and no-penetration along with isothermal conditions were enforced. 
Specifically, the boundary temperatures are set as $\Theta\left( t; x, 0\right) = \Theta\left( t; x, 1\right) = 0$.

The system of equations were solved numerically using finite differences on a staggered grid using the conservative form of the governing equations. 
A uniformly spaced staggered grid with $n_x \times n_y$ cells and spacing $\delta x = \delta y = h$ was adopted. 
The size of the computational cell is $\delta x = L_x/(L_y \times n_x)$ horizontally and $\delta y = 1/n_y$ vertically, where $n_x$ and $n_y$ are the number of horizontal and vertical grid points, respectively. 
A third-order Adam-Bashforth scheme was employed 
using a fixed time step $\delta t$. 
A second-order central differencing scheme was selected for the advection, diffusion and source terms. 
Pressure-projection based on the fractional step method was utilized to satisfy the incompressibility condition, which is manifested in the conservation of mass equation. 
Fast Fourier Transform (FFT) based algorithms were implemented to efficiently solve the parabolic and elliptic systems. 

\subsection{Continuous and Discrete Data Assimilation}

CDA is a nudging-based downscaling algorithm that employs coarse-scale observational data to constrain the large-scale features of the model solution in order to dynamically recover the fine-scale features.
This is achieved by introducing a source term, proportional to the discrepancy between the coarse scales of the model solution and those of the observational data, in the equations of motion and energy \citep{Azouani2013}.
CDA is applied here to the $2$D Rayleigh-B\'enard convection for downscaling coarse observational data by solving the following system of equations:
\begin{equation}
    \nabla \cdot \mathbf{w} = 0,
    \label{eqn:incompr_ContinuityCDA}
\end{equation}
\begin{equation}
    \frac{\partial \mathbf{w}}{\partial t} +  \left( \mathbf{w} \cdot \nabla \right) \mathbf{w} + \nabla p = \frac{Pr}{\sqrt{Ra}} \nabla ^2 \mathbf{w} + Pr \Psi \mathbf{e}_y + \mu_{\mathbf{u}} \left( I_{h^o} \left( \mathbf{u}^{o} \right) - I_{h^o} \left( \mathbf{w} \right) \right),
    \label{eqn:momentum_CDA}
\end{equation}
\begin{equation}
    \frac{\partial \Psi}{\partial t} + \left( \mathbf{w} \cdot \nabla \right) \Psi = \frac{1}{\sqrt{Ra}} \nabla ^2 \Psi + \mathbf{w} \cdot \mathbf{e}_y   + \mu_{\Theta} \left( I_{h^o} \left( \Theta^{o} \right) - I_{h^o} \left( \Psi \right) \right),
    \label{eqn:energy_CDA}
\end{equation}
where $\mathbf{w}=\left(\tilde u, \tilde v\right)$ is the downscaled velocity vector, $\Psi$ is the downscaled temperature, $\mathbf{\mu}=\left(\mu_{\mathbf{u}},\ \mu_{\Psi}\right)$ are non-negative constants called the nudging parameter, and $\mathbf{u}^{o}$ and $\Theta^{o}$ are the coarse-scale observations of the velocity and temperature, respectively. 
Note that the downscaling system is solved using the same boundary conditions as the reference system.
Finally, $I_{h^o}$ is an interpolation operator of an interpolant function $\phi$ for uniformly spaced measurements at a distance $h^o$ in $\Omega$, which is expressed as follows:
\begin{equation}
    I_{h^o}\left( \phi \left( \mathbf{x} \right) \right) = \sum_{k = 1}^{N_{h^o}} \phi \left( \mathbf{x}_k \right) \chi_{Q_k} \left( \mathbf{x} \right),
\end{equation}
where $\mathbf{x}=(x, y)$, $Q_k$ are disjoint subsets such that diam($Q_k$) $\leq h^o$ for $k=1,\ ...,\ N_{h^o}$, $N_{h^o}$ is the number of observational points, $\bigcup_{j=1}^{N_{h^o}} Q_j = \Omega$, $x_k \in Q_k$ and $\chi_E$ is the characteristic function of the set $E$.

It was proven that regardless of the initial conditions used for $\Psi$ and $\mathbf{w}$, the CDA solution converges exponentially to the true solution under the conditions \cite{Azouani2013} outlined.
Specifically, the system converges if sufficiently enough coarse observational data that is spaced at a sufficiently small $h^o$ are provided and the $\mathbf{\mu}$'s are appropriately selected.
CDA, thus, ensures that, for dissipative systems, the downscaled solution is synchronized with the reference solution for all times beyond the last observation time.

Inspired by the works of \cite{HAYDEN2011}, \cite{Azouani2013} and \cite{Celik2019}, we propose the DDA algorithm in which nudging is performed only at discrete instances when observations are available, while solving the original system at all other steps. 
Specifically, DDA simulates Equations (\ref{eqn:momentum_CDA}) and (\ref{eqn:energy_CDA}) when observations are available, and integrates the original Equations (\ref{eqn:momentum_true}) and (\ref{eqn:energy_true}) at all other times.
The discrete forms of the energy equations for CDA and DDA are presented in order to provide additional insight on the implementation of the proposed downscaling algorithms as: 
\begin{equation}
    \frac{ \Psi^{n+1} - \Psi^n }{\Delta t} + (  \mathbf{w} \cdot \nabla  \Psi )^n = \frac{1}{\sqrt{Ra}}  \nabla^2 \Psi^{n} + (\mathbf{w} \cdot \mathbf{e}_y)^{n} + \mu_{\Theta} \left[ I_{h^o}(\Theta^{o}(t^o_i)) - I_{h^o}(\Psi^n) \right]
    \label{eqn:CDA_Discrete}
\end{equation}
for CDA and 
\begin{eqnarray}
    \frac{ \Psi^{n+1} - \Psi^n }{\Delta t} + (  \mathbf{w} \cdot \nabla  \Psi )^n & = & \frac{1}{\sqrt{Ra}}  \nabla^2 \Psi^{n} + (\mathbf{w} \cdot \mathbf{e}_y)^{n} \nonumber \\
    & + &  \sum_{i=1}^{n_{t^o}} \delta\left(t - t_i^o\right) \mu_{\Theta} \left[ I_{h^o}(\Theta^{o}(t^o_i)) - I_{h^o}(\Psi^n) \right]
    \label{eqn:DDA_Discrete}
\end{eqnarray}
for DDA, where $\delta$ is the dirac delta function. In Equations (\ref{eqn:CDA_Discrete}) and (\ref{eqn:DDA_Discrete}) the algorithm advances from time step $n$ to $n+1$ and $t^o_i$ is the instance of the $i^{th}$ observation.
CDA, therefore, applies the nudging term that employs the $i^{th}$ observation between $t_i^o$ and $t_{i+1}^o$.
On the other hand, DDA includes the nudging term only at the time step when the observation was made and solves the original system of equations otherwise.  

The assimilation of perfect observations was analyzed both theoretically and numerically for CDA. 
The effect of noisy observations on the convergence rate and the asymptotic error level have only been examined theoretically by \cite{Bessaih2015}.
\cite{Bessaih2015} showed that the expected $\ell_2$-norm of the error, $\Lambda$, between the downscaled ($Z_d$) and true ($Z_t$) solution fields is proportional to the trace of the covariance of the observational noise, which is expressed as follows:
\begin{equation}
    \Lambda \coloneqq  \mathbb{E} \left[ \int \left( Z_{d} - Z_{t} \right)^2 d \Omega \right] \propto C \cdot Tr\left[ Cov \right]
    \label{eqn:measure}
\end{equation}
where $Z$ is an arbitrary variable, $d \Omega$ represents a differential area in the domain, $C$ is a functional that depends on the numerical and downscaling parameters, and $Cov$ is the covariance of the noise. 
One objective of this work is to numerically investigate the asymptotic error level and the convergence rate. 
Another objective is to estimate the functional dependence of $\Lambda$ on the covariance of the noise, and spatial and temporal resolutions of the observations, based on which the CDA and DDA algorithms were compared.

\subsection{Generation of Random Observations}
\label{sec:randomObs}

The solution of Equations (\ref{eqn:incompr_Continuity})-(\ref{eqn:energy_true}) with the initial and boundary conditions described in Section \ref{sec:rb2d} is considered as the reference solution.  Coarse observations of the system states were generated by subsampling points every $h^o$ spatially and $\delta t^o$ temporally from the reference solution.  (Note that as further discussed below, the reference solution is obtained on a sufficiently fine computational grid that it can be practically treated
as the exact solution.) The spatial and temporal resolutions of the observational grid are characterized using the non-dimensional ratios, called spatial and temporal downscaling factors, defined as $\mathcal{R} \equiv h^o / h$ and $\mathcal{S} \equiv \delta t^o / \delta t$, respectively.

The coarse-scale observations of the reference solution are then perturbed by adding randomly sampled independent and identically distributed (iid) white noise sampled from a Gaussian distribution with zero mean and standard deviation $\sigma$.
In this work, $\bm{\sigma}=\left(\sigma_T, \sigma_{\mathbf{u}} \right)$ is also referred to as the noise level, where the noise of the two variables are independent to mimic different measurement sources, a practical scenario encountered with observational data.
Each noisy observed realization is then downscaled separately to generate an ensemble of downscaled solutions, each representing a possible solution trajectory of the system. 
The performance of the downscaling algorithms was then evaluated quantitatively based on the asymptotic error level between the downscaled and reference solutions, and the rate at which they reach the asymptotic behavior, i.e.\ the convergence rate.

\section{Experimental Setup}
\label{sec:expSet}

The experiments were conducted using the Rayleigh-B\'enard convection with a high Rayleigh number of $Ra=2 \times 10^8$ and $Pr=0.7$, which corresponds to a highly chaotic system that mimics environmental flows \citep{Weidauer_2010, Chill2012}.
Numerically, the system of equations are solved on a uniformly spaced staggered grid composed of $1200 \times 400$ grid points (for scalar fields) discretizing a domain of length $L_x = 3$ and height $L_y = 1$. 
The time step $\delta t$ was set to $5 \times 10^{-4}$, leading to a maximum Courant number of about $0.15$ and a maximum grid Reynold number of 
$\mathcal{O}(10)$.  This ensures that fine scale dynamics are suitably captured.  In other words, the computed solution is essentially insensitive to 
further refinement of the spatial or temporal resolution of the computational grid.

The sensitivity of the convergence rate and the asymptotic error level, of both CDA and DDA, to the downscaling parameters and the observation noise level were investigated to determine the corresponding functional dependencies. 
Specifically, the impact of the observation noise level was numerically studied by analyzing the sensitivity of the downscaled solution statistics to different values of $\bm{\sigma}$.
The solution variables were initialized with a random field and the computations were carried out over a sufficiently long interval to observe convergence. 
This generally occurred within $t = 25$ with DDA, and for a much shorter window ($t = 3$) with CDA. 
Moreover, the asymptotic error level was analyzed for the final steps, which, in this study, correspond to $t_f^{DDA} = 49.9$ for the DDA and $t_f^{CDA} = 15$ for the CDA.

A $0^{th}$ order interpolant operator was selected because it is computationally cheapest with no noticeable impact on the convergence rate, as was also suggested by \citep{Altaf2017}. 
The nudging factor for both the energy and momentum equations were chosen the same with $\mu_{\Theta}=\mu_{\mathbf{u}} = \mu$ to limit the number of hyper-parameters.  In the present work, a grid search was performed to find the best $\mu$, where the downscaled fields were ensured to converge to the reference solution. Note that the theory in~\citep{Azouani2013} predicts the existence of a suitable range of parameters, and computational experience~\citep{Altaf2017} suggests that CDA predictions always converge so long as $\mu$ is selected within this range.  The 
performance of the downscaling algorithms also depends on the spatial and temporal resolutions of the available observations.
The sensitivity of both CDA and DDA to these parameters were, therefore, also investigated.

The quality of the downscaled solutions were assessed by computing various skill scores to quantify the errors relative to the reference solution and variability relative to the mean solution.
The skill scores were computed for each of the solution variables and included the root-mean-squared-error, which is expressed as:
\begin{equation}
    RMSE(t) = \sqrt{ \frac{1}{n_x n_y} \Sigma_{i=1}^{n_x  n_y} \left[ Z_d (\bm{x}_i; t) - Z_{t}(\bm{x}_i; t) \right]^2},
\end{equation}
where $\bm{x}_i = (x_i, y_i)$ is the position of the grid point and $t$ is time; the relative root-mean-squared error:
\begin{equation}
    RRMSE(t) = \frac{\sqrt{ \Sigma_{i=1}^{n_x  n_y} \left[ Z_{d} (\bm{x}_i; t) - Z_{t}(\bm{x}_i; t) \right]^2}}{\sqrt{ \Sigma_{i=1}^{n_{x} n_{y}} \left[ Z_{true}(\bm{x}_i; t) \right]^2}}.
\end{equation}
and the absolute error:
\begin{equation}
    AE(t) = \frac{1}{n_x n_y} \Sigma_{i=1}^{n_{x}  n_{y}}  \vert Z_{d} (\bm{x}_i; t) - Z_{t}(\bm{x}_i; t) \vert.
\end{equation}
The variability of the downscaled solutions was evaluated as the average ensemble spread as follows:
\begin{equation}
    AES(t) = \sqrt{\frac{1}{N_{e} -1} \Sigma_{k=1}^{N_{e}} \Sigma_{i=1}^{n_{x}n_{y}} \left(  Z_{d}^k (\bm{x}_i; t) - \overline{ Z_{d} }(\bm{x}_i; t) \right)^2 },
\end{equation}
where $N_{e}$ is the number realizations forming the ensemble of downscaled solutions and $\overline{Z_{d}}(\bm{x}_i; t)$ is the average of the samples at $\bm{x}_i$ and time $t$. 
Finally, the measure introduced by \cite{Bessaih2015} on the expected $\ell_2$-error of the downscaled members, denoted by $\Lambda$, was also computed based on the definition in Equation (\ref{eqn:measure}) for the last time step of the simulation.
This enables us to numerically assess the theoretical results and to further investigate the functional dependence on downscaling parameters.

\noindent {\bf Remark:} It is instructive to contrast the resolution of the coarse observation grids to those of the integral time scales of the flow.  Based on 
the definitions above, the spatial and temporal resolutions of the observations grid are given by $h^0 = {\cal R} h$ and $\delta t^0 = {\cal S} \delta t$.  In the
downscaling experiments conducted below, the downscaling ratios are selected such that $ 5 \leq {\cal R} \leq 25$ whereas $ 1\leq {\cal S} \leq 50$.  With these choices, $h^0$ and $\delta t^0$ remain substantially smaller than the flow integral space scale and time scale, respectively.  In particular, the 
dimensionless large-scale turnover time, $\tau \equiv 2L_y/u_{max}$, estimated based on the reference solution, is approximately 3, 
which indicates that $\delta t^0 / \tau$ remains much smaller than 1 for the entire range of the ${\cal S}$ values considered.

\section{Results}
\label{sec:results}

This section presents and discusses the results of the numerical experiments illustrating the implementation of the DDA and CDA algorithms to downscale coarse scale noisy observations of the Rayleigh-B\'enard flows. 
The DDA results are mostly presented in the main text and when not included the CDA results are mirrored in the supplementary material, but are analyzed in the text for direct comparison.

\subsection{Sensitivity to Noise Levels}
\label{res:noiseLevels}

The sensitivity of each of the CDA and DDA algorithms to observation noise levels were examined by first selecting a nominal case, defined by a nominal noise level equal to $10 \%$ of the range of the solution variable at the final time step. The noise levels for the nominal case were set to $(\sigma_T, \sigma_{\mathbf{u}}) = (0.1, 0.05)$. Higher and lower noise levels are then considered by scaling both noise levels simultaneously. Consequently, the results presented indicate $\sigma_T$ only.

Noisy realizations of the coarse-scale observations of the system are sampled as described in Section \ref{sec:randomObs} and are then downscaled using CDA and DDA, respectively.
The DDA and CDA nudging parameters are set to $7.0$ and $3.0$, respectively, for an observational grid resolution of $\mathcal{S} = 10$ and $\mathcal{R} = 5$. 
The nudging factor was tuned by solving the downscaling equations for multiple values of $\mu$, and the solution with lowest RRMSE at the final time step is used throughout.
Downscaling the noisy realizations was conducted starting from a random field to verify that, given sufficient observations and appropriate nudging, both algorithms converge regardless of the initial conditions.

\begin{figure}[!htbp]
    \centering
    \subfloat[AE]{\includegraphics[width=60mm]{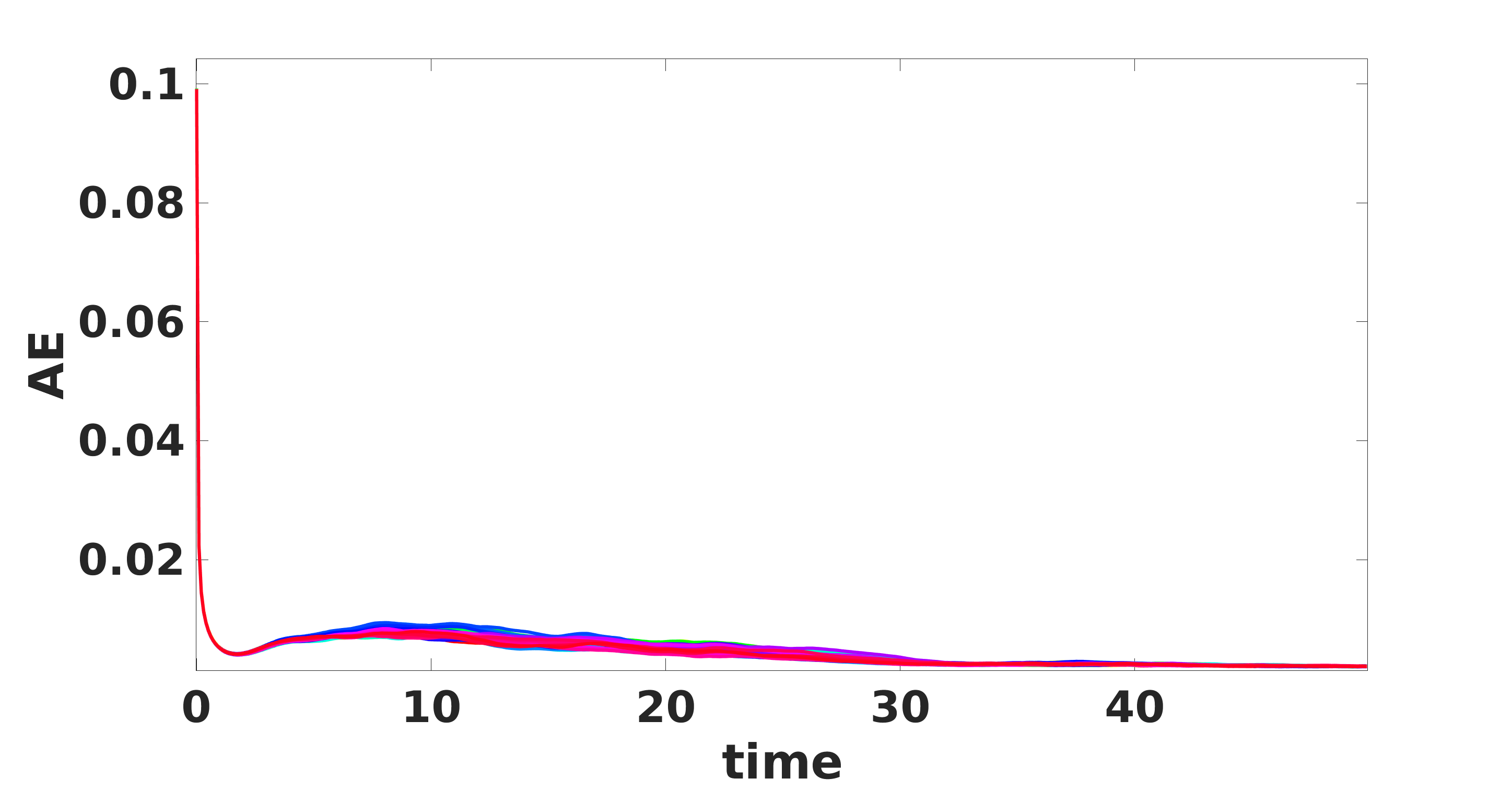}\label{part:aFIG1}} \,
    \subfloat[RMSE]{\includegraphics[width=60mm]{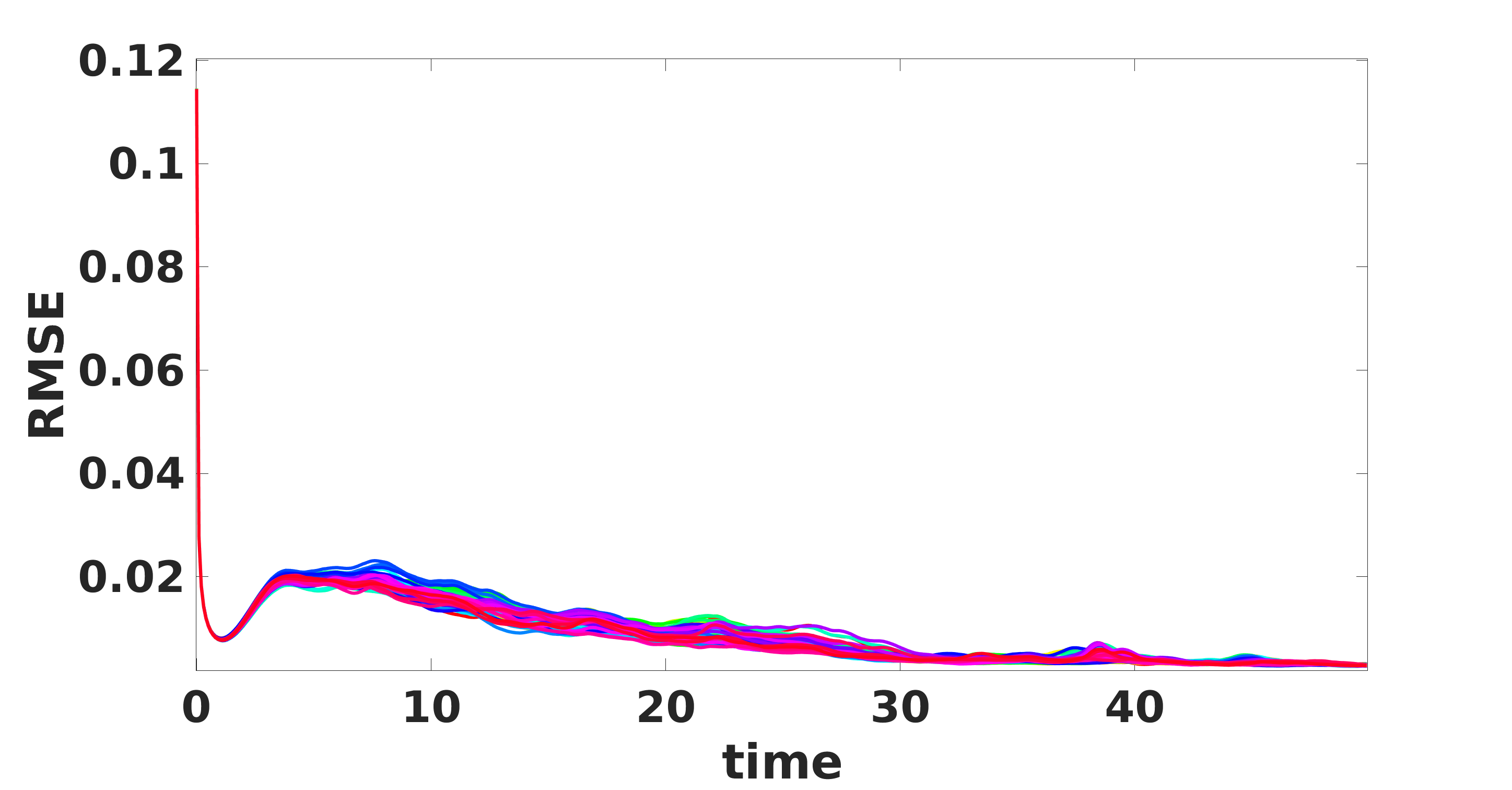}\label{part:bFIG1}} \,
    \subfloat[AES]{\includegraphics[width=60mm]{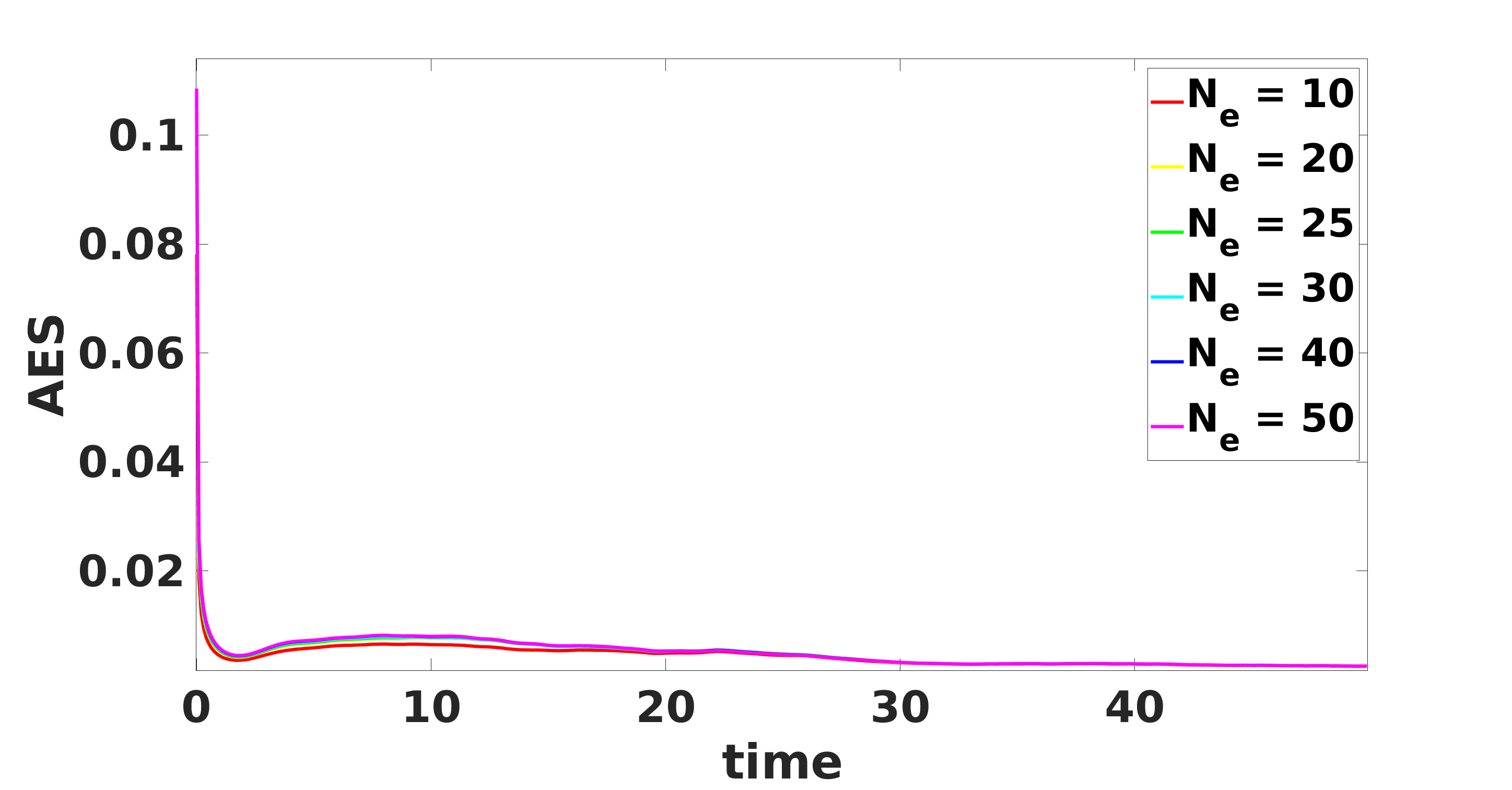}\label{part:cFIG1}} \,
    \subfloat[RRMSE - DDA]{\includegraphics[width=60mm]{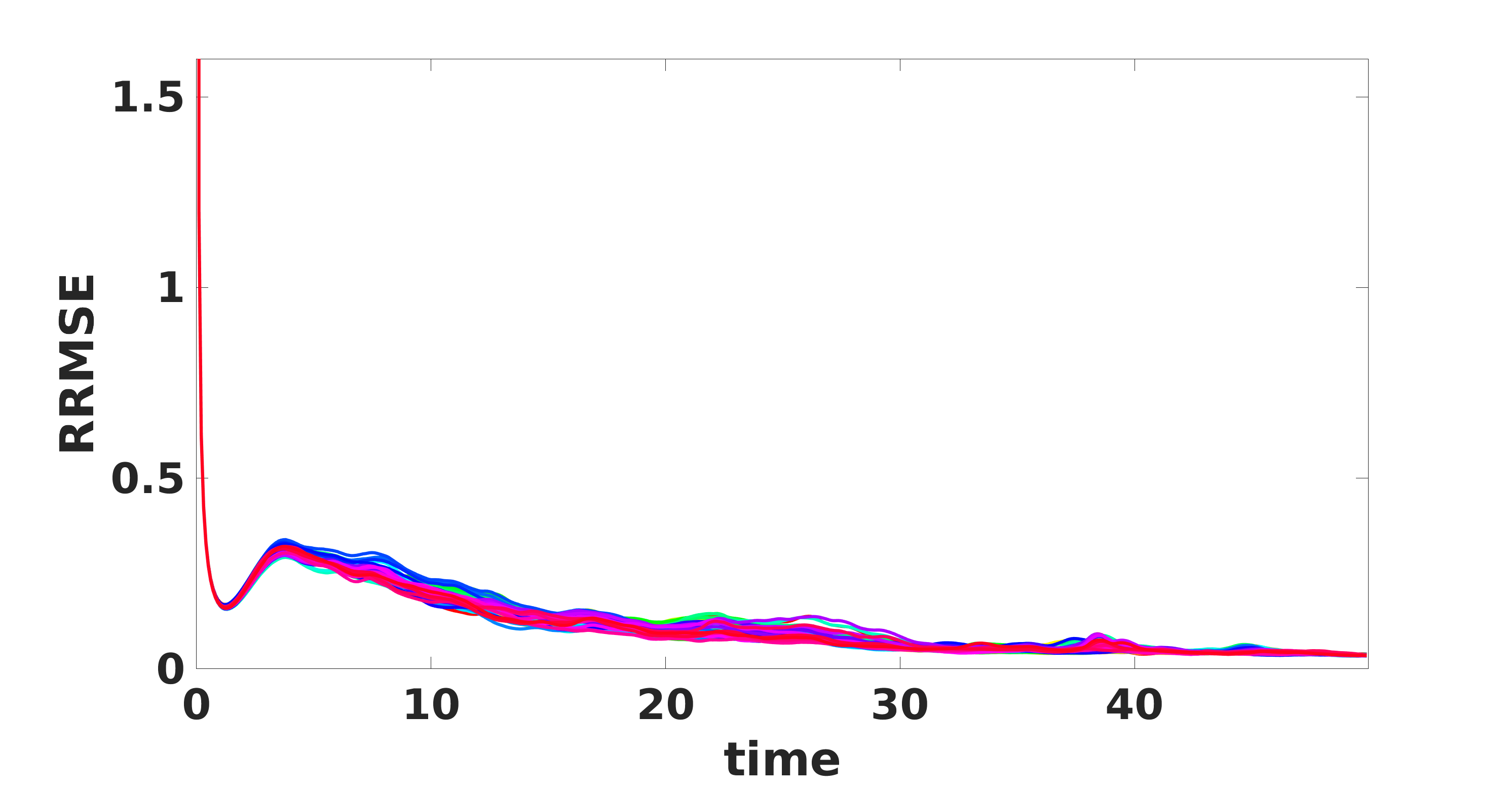}\label{part:eFIG1}} \,
    \subfloat[RRMSE - CDA]{\includegraphics[width=60mm]{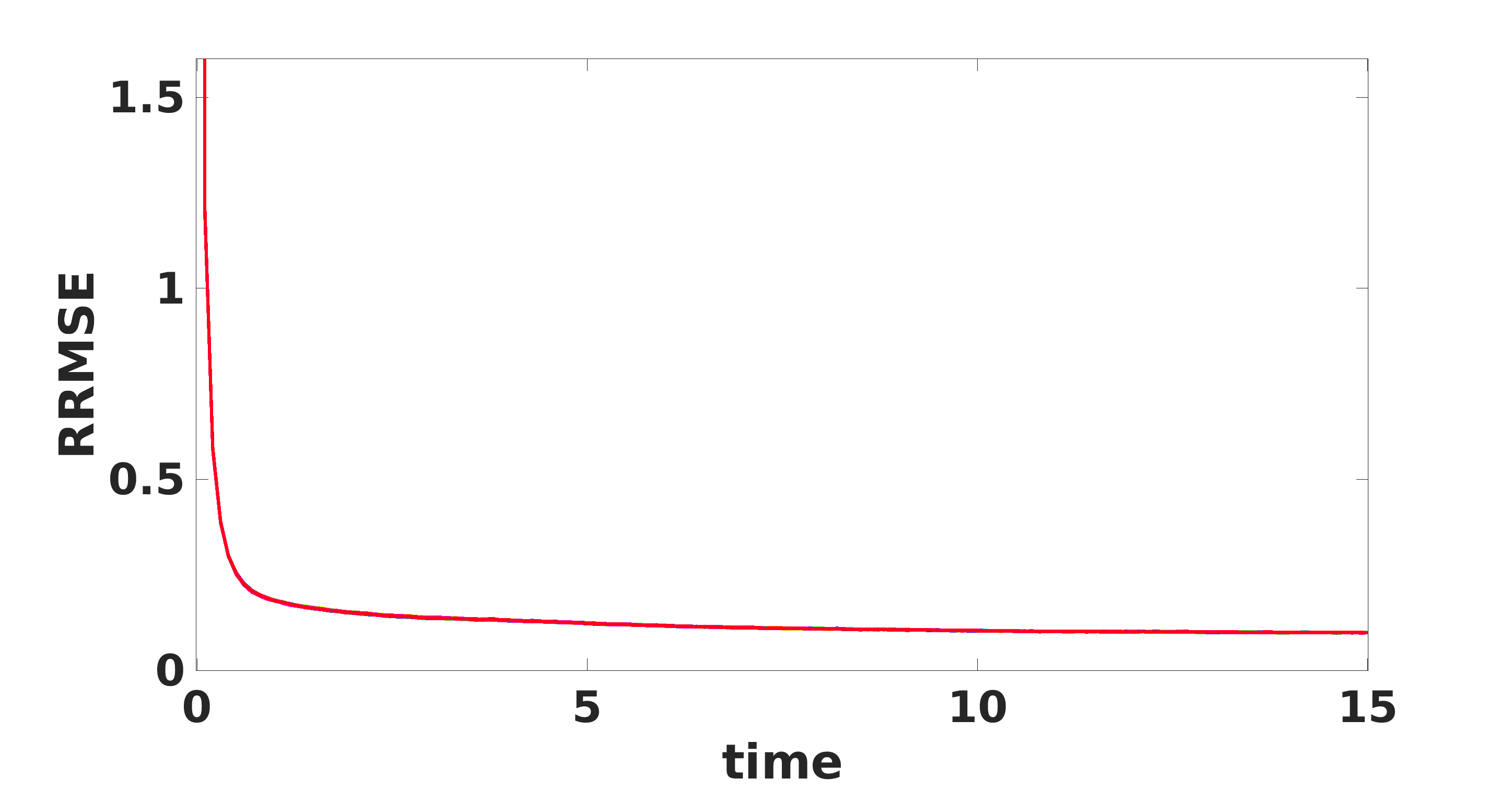}\label{part:fFIG1}} \,

\caption{Evolution of different skill scores corresponding to the DDA downscaled temperature fields, each member is driven by observations of the nominal noise level of $\left( \sigma_T, \sigma_{\mathbf{u}}\right) = \left( 0.1, 0.05\right)$. Each plot contains curves corresponding to $50$ downscaled samples, the AES is presented for different number of realizations as indicated.}
\label{Fig:nominalNoiseLvl}
\end{figure}

Figure \ref{Fig:nominalNoiseLvl} presents the time evolution of the AE, RMSE, RRMSE of the DDA downscaled temperature solutions.
The figure also illustrates the AES for ensembles composed of $10$, $20$, $25$, $30$, $40$ and $50$ downscaled solutions. 
The evolution of the RRMSE from the CDA simulation is also shown for comparison. 
The results of the velocity components are similar to those of the temperature; they are briefly discussed in this section with the corresponding results 
provided in the Supplementary.

The time evolution of both the AE and RMSE of the DDA downscaled solutions for all $50$ members of the ensemble are illustrated in Figures \ref{part:aFIG1} and \ref{part:bFIG1}, respectively. 
The skill scores decrease in time until reaching a plateau at approximately $t=30$. 
Figure \ref{part:cFIG1} shows a rapid decay of the AES, for all considered ensemble sizes, following similar patterns to those of AE and RMSE. 
The plots also indicate that the smallest AES values are at the plateau. 
Similar to the downscaled temperature fields, the errors of the downscaled velocity fields also decay to an asymptotic error level, and achieve AES values that are comparable to those of the downscaled temperature fields.  

Figures \ref{part:eFIG1} and \ref{part:fFIG1} illustrate the decay of the RRMSE for both DDA and CDA downscaled fields. 
These plots show that with CDA the RRMSE decays smoothly and exponentially in time, whereas with DDA the RRMSE decays non-monotonically and at a slower rate. 
As time evolves, the skill scores of both algorithms asymptote to a plateau, with the asymptotic error skill scores corresponding to DDA being lower than that of CDA. 
Specifically, the asymptotic RRMSE of the DDA downscaled fields were approximately three times smaller than that corresponding to CDA, and is approximately one third of the observations' noise level.
The RRMSE of the downscaled velocity components experience a similar evolution to that of the downscaled temperature fields, where the RRMSE decays smoothly and exponentially with CDA and non-monotonically with DDA.
This error behavior can be attributed to the nature of the nudging involved with each algorithm. 
Specifically, CDA continuously nudges the system of equations, which introduces observational noise at every integration time step.
On the other hand, DDA nudges the system of equations at discrete times, allowing the unforced system of equations to relax this sudden change imposed by the nudging at discrete time steps.
Hereafter, we only discuss the evolution of the RRMSE of the downscaled temperature solution because the skill scores evolve similarly in time.

\begin{figure}[!htbp]
    \centering
    \subfloat[$\mathbf{t = 0}$]{\includegraphics[width=120mm]{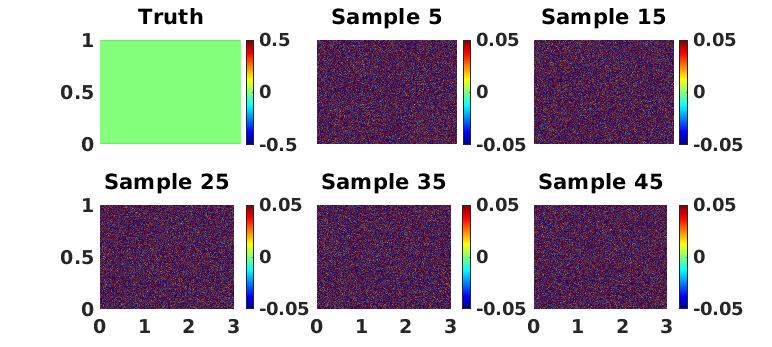}} \,
    \subfloat[$\mathbf{t = 10}$]{\includegraphics[width=120mm]{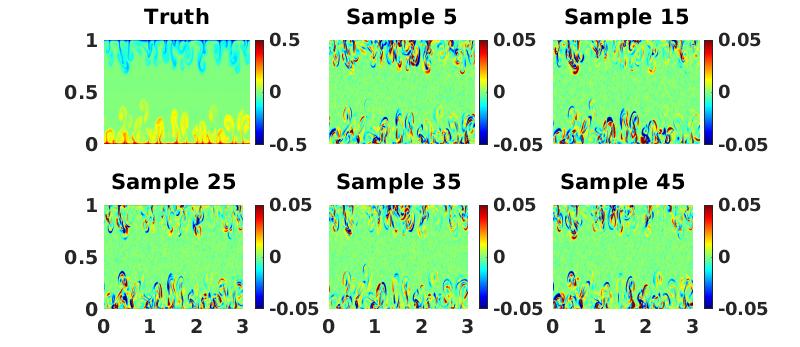}} \,
    \subfloat[$\mathbf{t = 45}$]{\includegraphics[width=120mm]{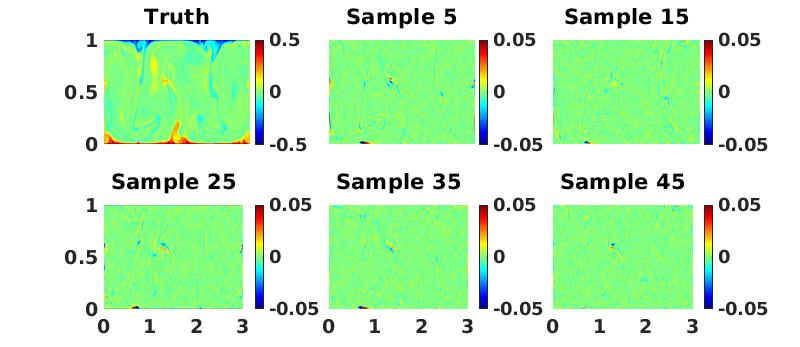}} \,

\caption{Instantaneous snapshots of the temperature distribution for the noise-free Rayleigh-B\'enard simulation along with the  difference between the downscaled sample solutions and the truth. Snapshots are shown for the initial time step (top), an intermediate time step (middle) and a time step when the algorithm had converged (bottom).}
\label{Fig:nominalTempProfile}
\end{figure}

The reference temperature solution along with its difference with the DDA downscaled solutions are presented at different time instants in Figure \ref{Fig:nominalTempProfile}. 
The snapshots first show the initial condition of the reference system and the random initialization of the DDA systems, which is expressed by the noisy error field.
The solution at the intermediate time shows a decrease in the errors across the domain, mostly concentrated near the plumes, and are of the order of the observation noise level.
At a later time, the errors practically vanish across the domain, with the highest error values being an order of magnitude smaller than the noise level. 
The highest errors occur near detachment and mixing zones due to the more complex fine-scale dynamics, as expected.

\begin{figure}[!htbp]
    \centering
    \subfloat[$\mathbf{t = 15}$]{\includegraphics[width=120mm]{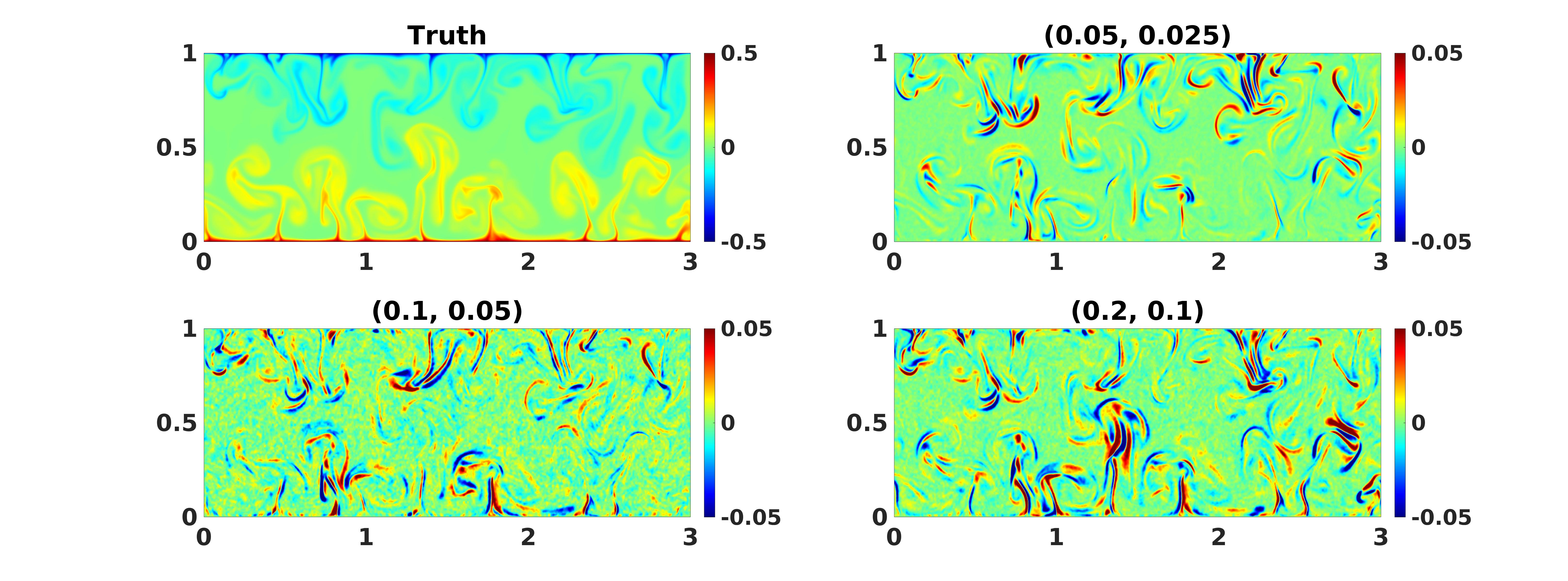}} \,
    \subfloat[$\mathbf{t = 45}$]{\includegraphics[width=120mm]{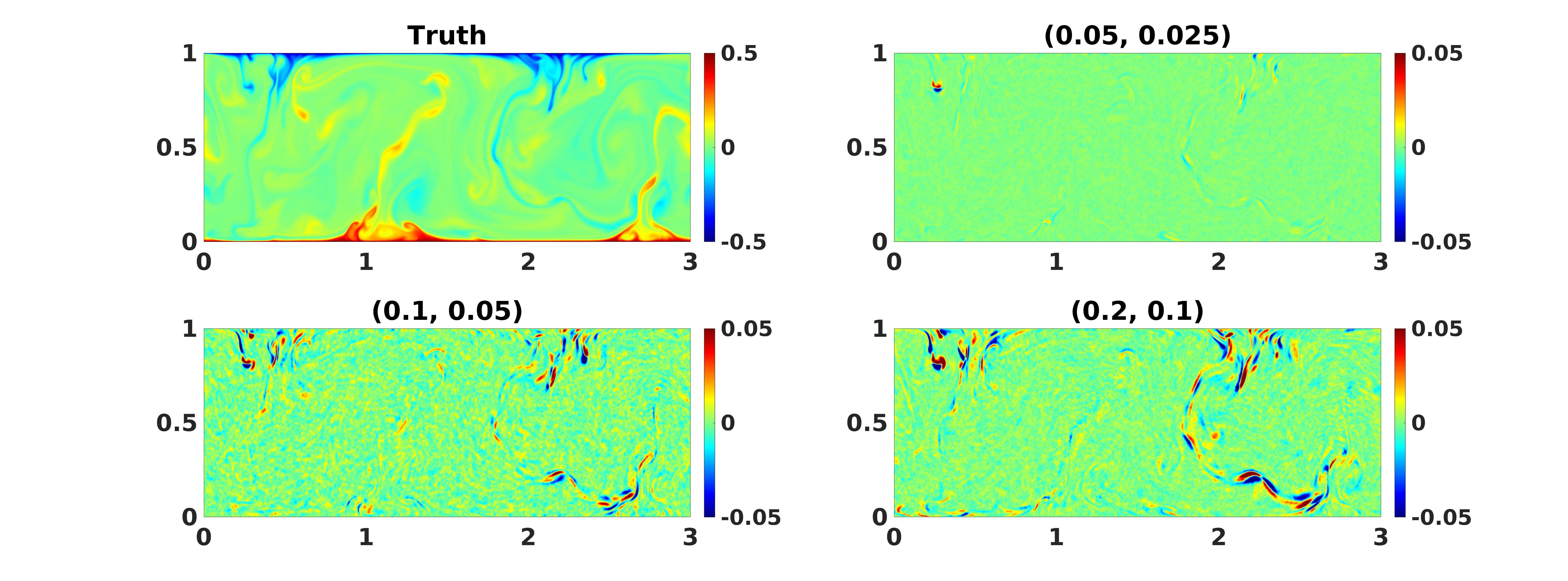}} \,

\caption{Instantaneous snapshots of the temperature distribution for the noise-free Rayleigh-B\'enard simulation along with the difference between the downscaled solutions and the truth.  Snapshots are shown for an intermediate time step (top) and a time step when the algorithm had converged (bottom).
At both instants, plots are generated for different noise levels, as indicated.}
\label{Fig:TempProfiles_diffNoiseLvls}
\end{figure}

To provide insight into the effect of the noise level on the downscaled fields, the reference temperature field and the difference of the DDA downscaled temperature fields with this reference are presented for different snapshots in Figure \ref{Fig:TempProfiles_diffNoiseLvls}.  In all cases,
the results indicate that the difference between the downscaled and reference temperature fields is largest around plumes and 
their corresponding boundaries.   At an intermediate time ($t=15$), the discrepancy fields reveal structures whose shapes and amplitudes depend
on the observation noise level, especially when the latter is large.  This is not surprising because in this situation the downscaled solution is still 
far from the reference solution.  
As further observations are assimilated in time ($t=45$), the impact of increasing the noise level essentially manifests itself as an amplitude 
scaling of the discrepancy field.  Thus, when the error plateau is reached, the discrepancy fields exhibit similar structure, with amplitude differences that
scale with the observation noise level.

\begin{figure}[!htbp]
    \centering
    \subfloat[Absolute Error]{\includegraphics[width=74mm]{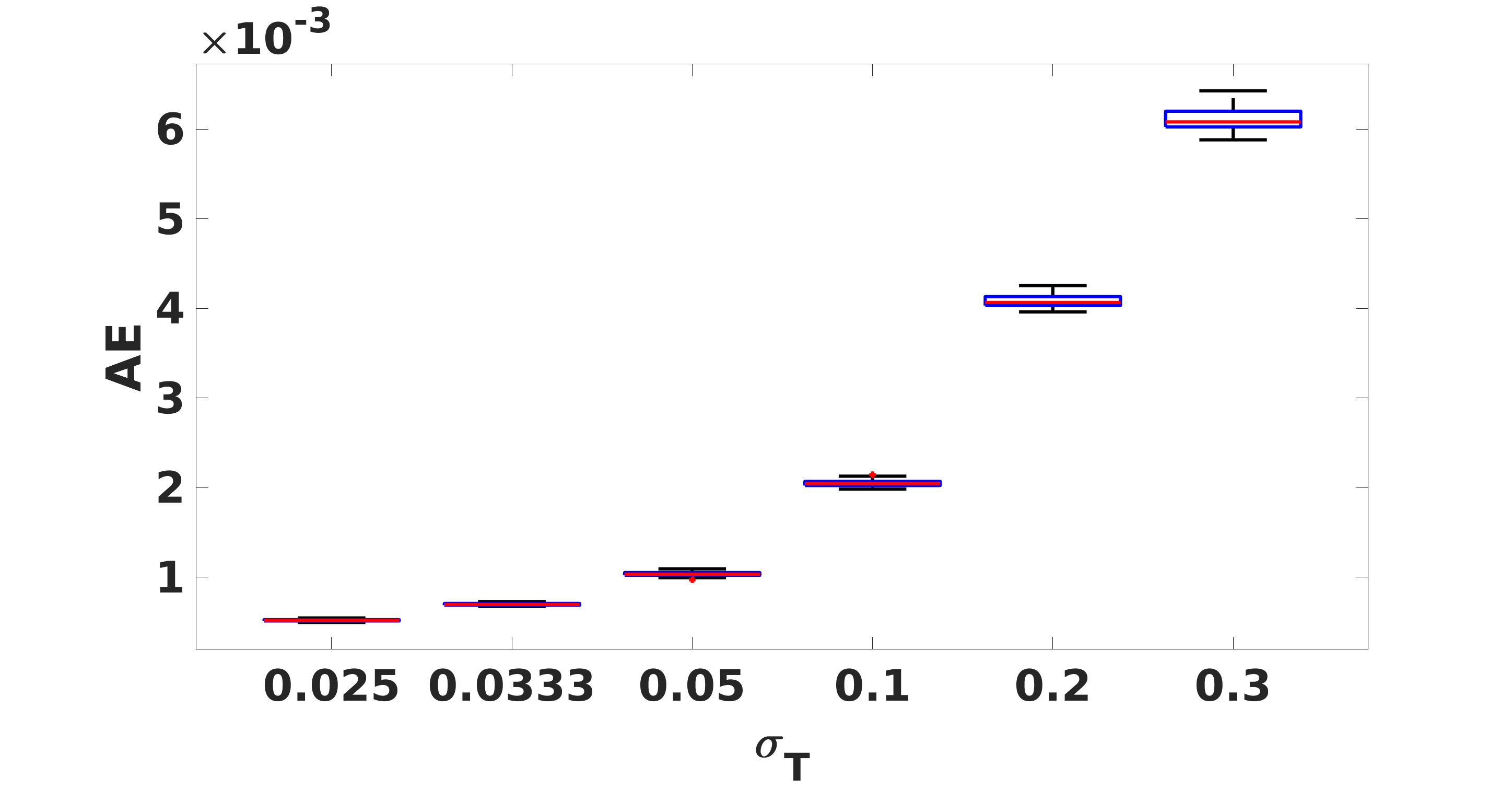}\label{part:FIG3a}} \,
    \subfloat[RRMSE]{\includegraphics[width=74mm]{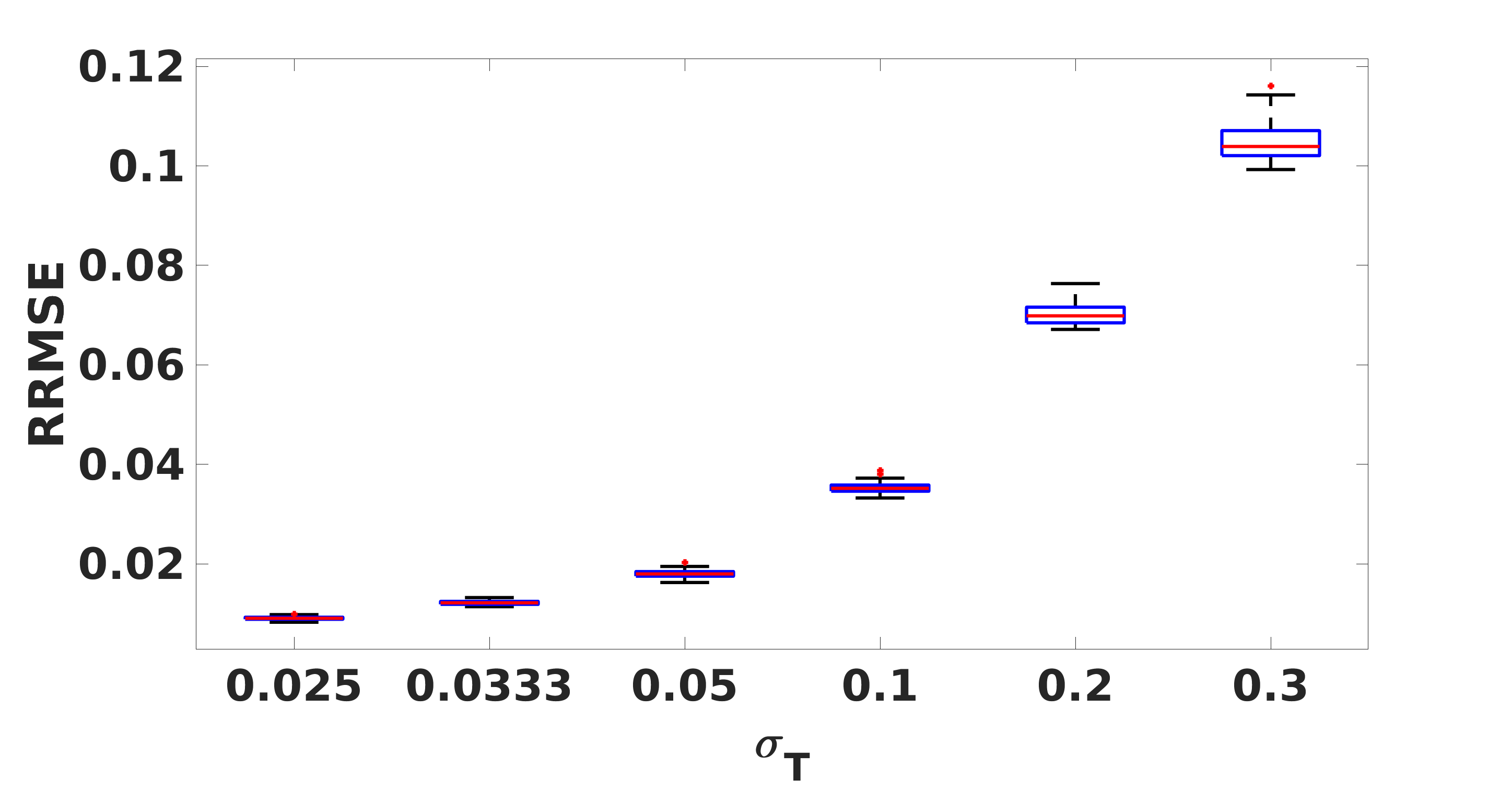}\label{part:FIG3b}} \,
    \subfloat[AES]{\includegraphics[width=74mm]{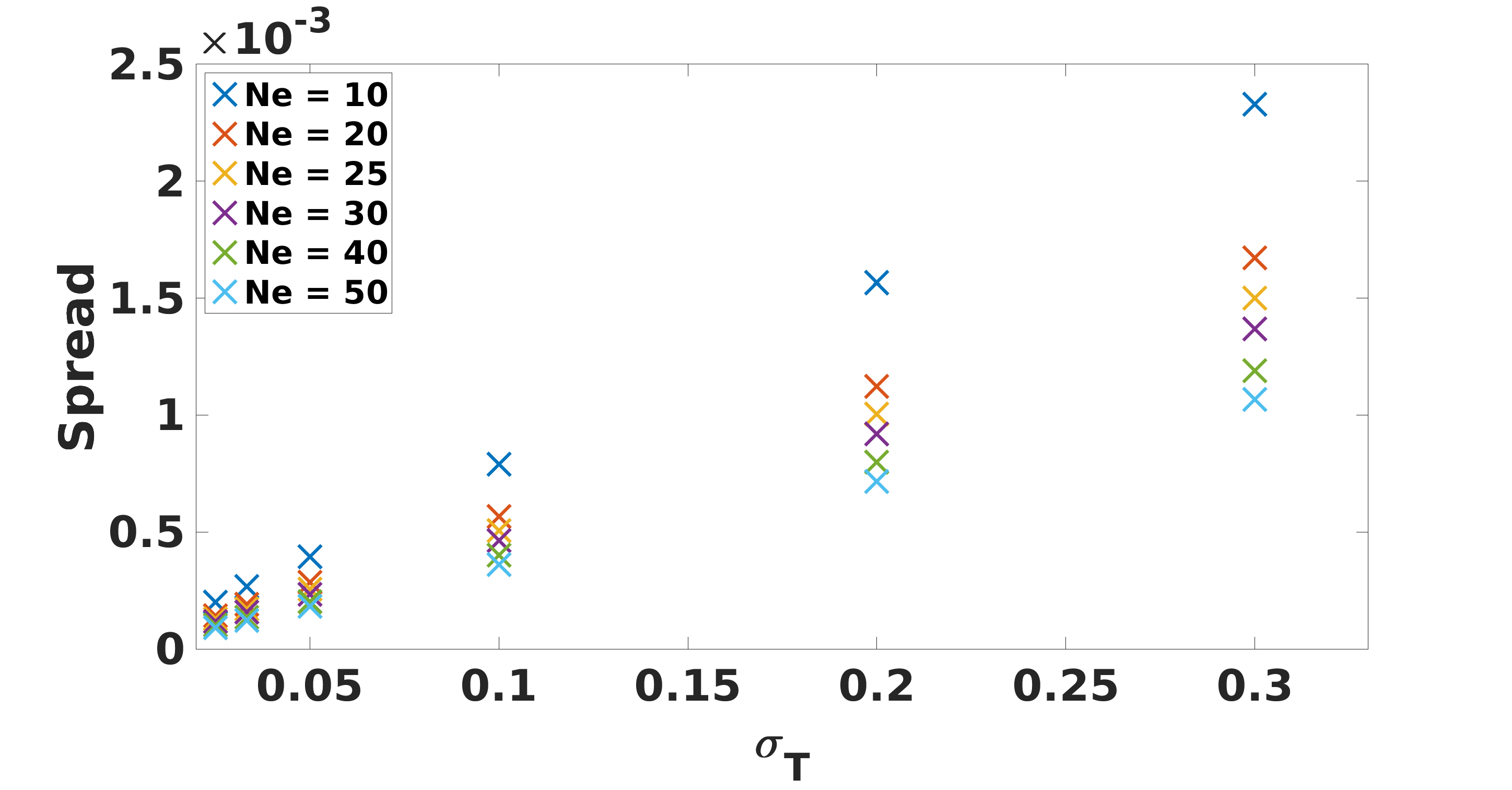}\label{part:FIG3c}} \,
    \subfloat[$\Lambda$]{\includegraphics[width=74mm]{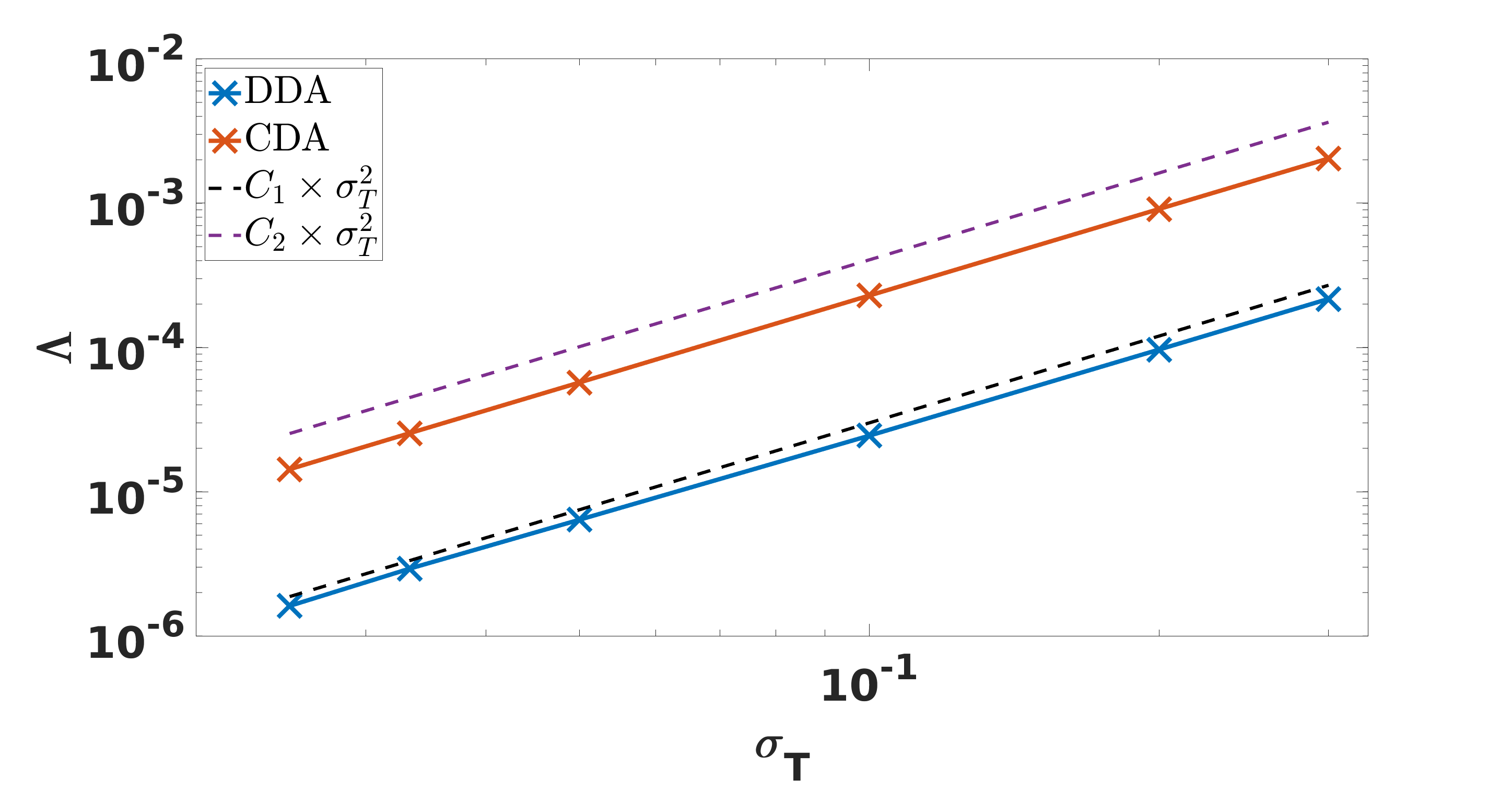}\label{part:FIG3d}} \,

\caption{Distribution of the skill scores (indicated) based on the $50$ sample solutions for the ensemble of temperature solutions of the (a)-(d) DDA solution, and (d) CDA solution at their respective final time steps. Sample solutions were generated by varying the noise levels (indicated by the horizontal axis), while keeping the remaining parameters fixed.}
\label{lastStep_noiseLevels}
\end{figure}

The distributions of the asymptotic skill scores corresponding to the DDA downscaled temperature field at the final time step are shown in Figure \ref{lastStep_noiseLevels}. 
The box plots in Figures \ref{part:FIG3a} and \ref{part:FIG3b} depict an increase in the asymptotic AE and RMSE quadratically with the noise level. 
Figure \ref{part:FIG3c} shows the decrease in the AES for increasing ensemble size at a rate of approximately $1/\sqrt{N_{ens}}$, as expected for random sampling of the observation noise. The figure also illustrates a linear decrease of the AES with decreasing noise level and a fixed ensemble size. 
Finally, the curves in Figure \ref{part:FIG3d} present the increase of $\Lambda$ with $\sigma$ on a log plot for both CDA and DDA.
The results suggest that $\Lambda$ is proportional to the variance of the noise for both the CDA and the DDA downscaled fields, consistent with the theory of \cite{Bessaih2015}.

\subsection{Sensitivity to the Observation's Temporal Sampling Rate}
\label{res:freqObs}

This section analyzes the impact of the temporal observation sampling rate on the DDA and CDA solutions.
Numerical experiments were performed with different observation frequencies by systematically varying $\mathcal{S}$ within the broad range $1 \leq \mathcal{S} \leq 50$. 
In doing so, the noise level was held fixed with $\bm{\sigma} = (0.1,0.05)$ and $\mathcal{R} = 5$.

\begin{figure}[!htbp]
    \centering
    \subfloat[$\mathcal{S} = 1$]{\includegraphics[width=110mm]{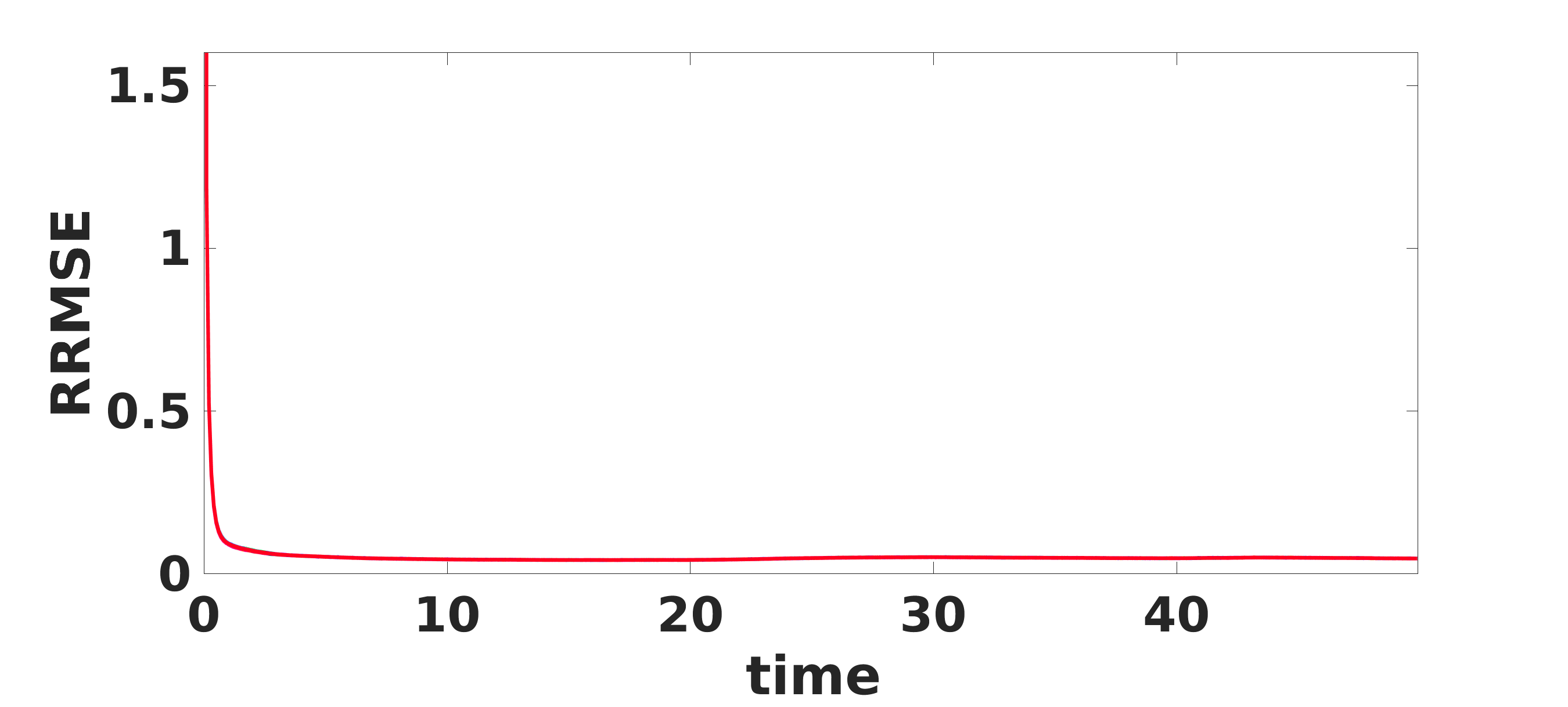}} \,
    \subfloat[$\mathcal{S} = 10$]{\includegraphics[width=110mm]{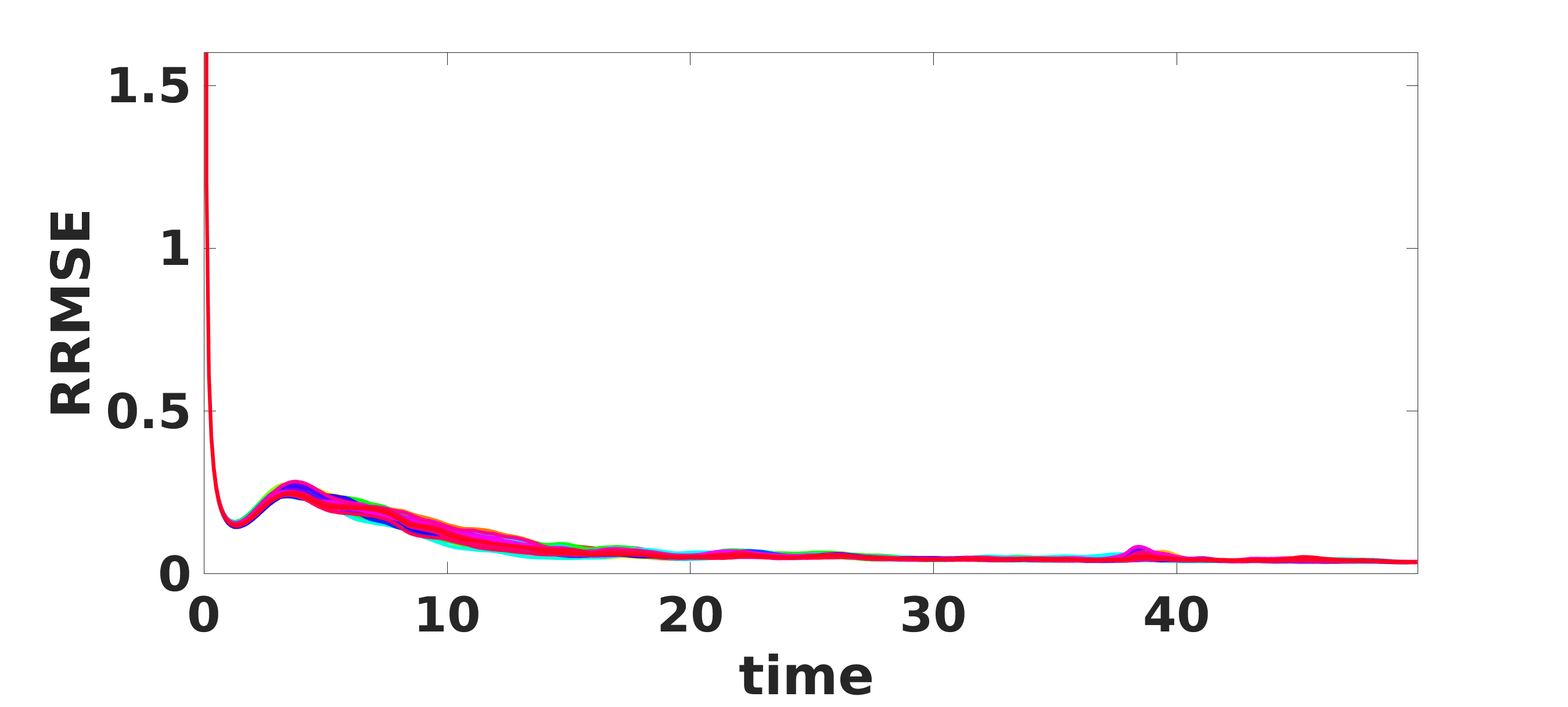}} \,
    \subfloat[$\mathcal{S} = 24$]{\includegraphics[width=110mm]{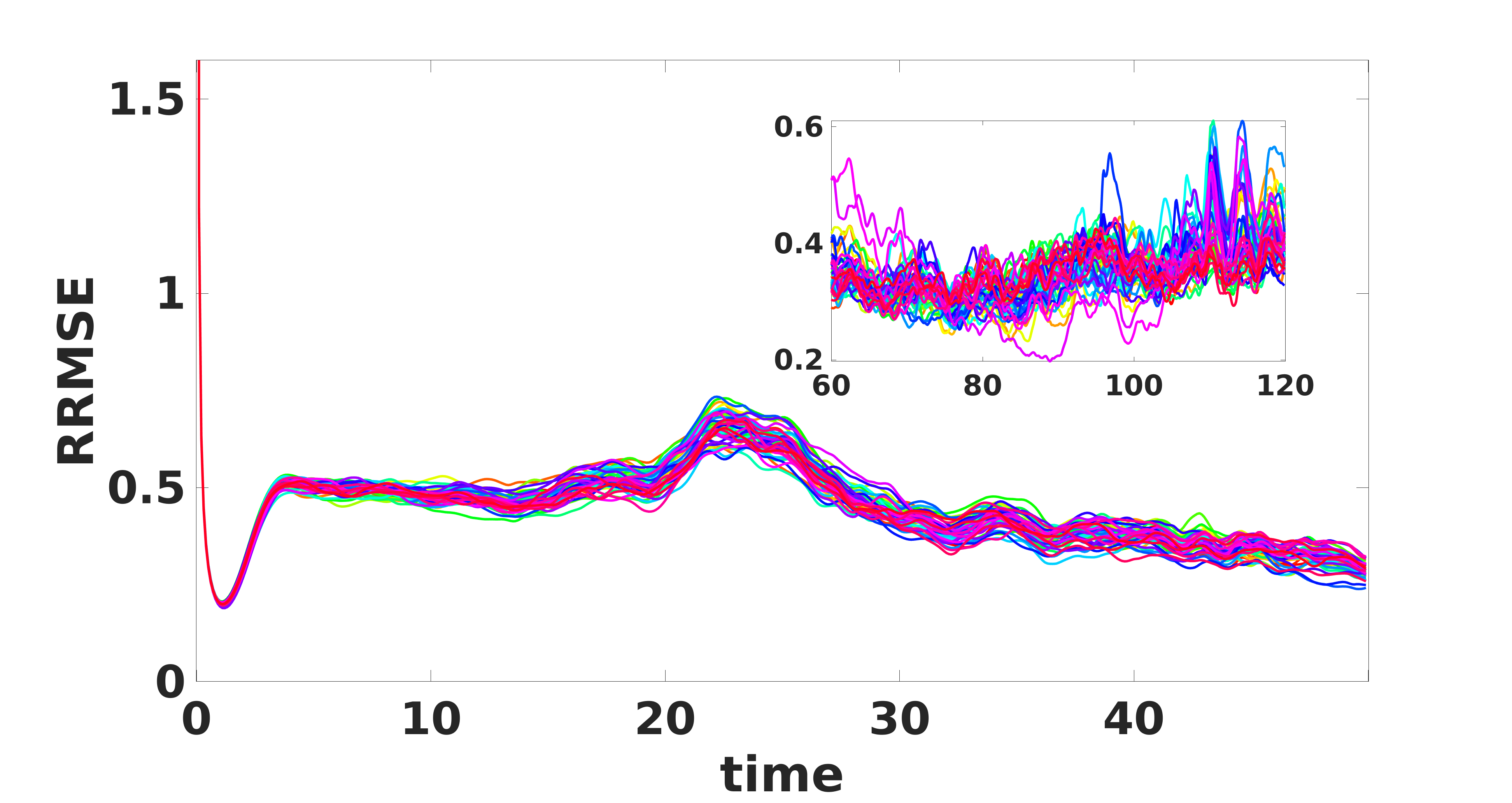}} \,

\caption{Evolution of the RRMSE of $50$ individually downscaled temperature solutions for: (a) $\mathcal{S} = 1$, (b) $\mathcal{S} = 10$ and (c) $\mathcal{S} = 24$.  Note that with $\mathcal{S}=24$ the simulations were carried out till $t=120$, as reflected in the inset.}

\label{evo_RRMSE_freqObs}
\end{figure}

Figure \ref{evo_RRMSE_freqObs} shows the time evolution of the temperature RRMSE resulting from DDA for $\mathcal{S}=1$, $10$ and $24$. 
For $\mathcal{S}=1$, the RRMSEs decay monotonically to a plateau. 
For $\mathcal{S}=10$ and $24$, however, the behavior becomes non-monotonic. 
In both cases, the RRMSE initially drops rapidly, then increases before decreasing again to eventually reach a plateau. 
As $\mathcal{S}$ increases, the downscaled fields exhibit a wider spread. 
Specifically, for $\mathcal{S} = 1$, the RRMSE curves are almost indistinguishable, whereas for larger $\mathcal{S}$, the variability in the downscaled fields becomes substantial, though a smaller mean error is achieved for $\mathcal{S}=10$.
CDA also seems to be less sensitive to $\mathcal{S}$ compared to DDA, where the time evolution of the temperature RRMSE of different members is almost identical and the plateau is reached almost immediately.

\begin{figure}[!htbp]
    \centering
    \subfloat[Absolute Error]{\includegraphics[width=74mm]{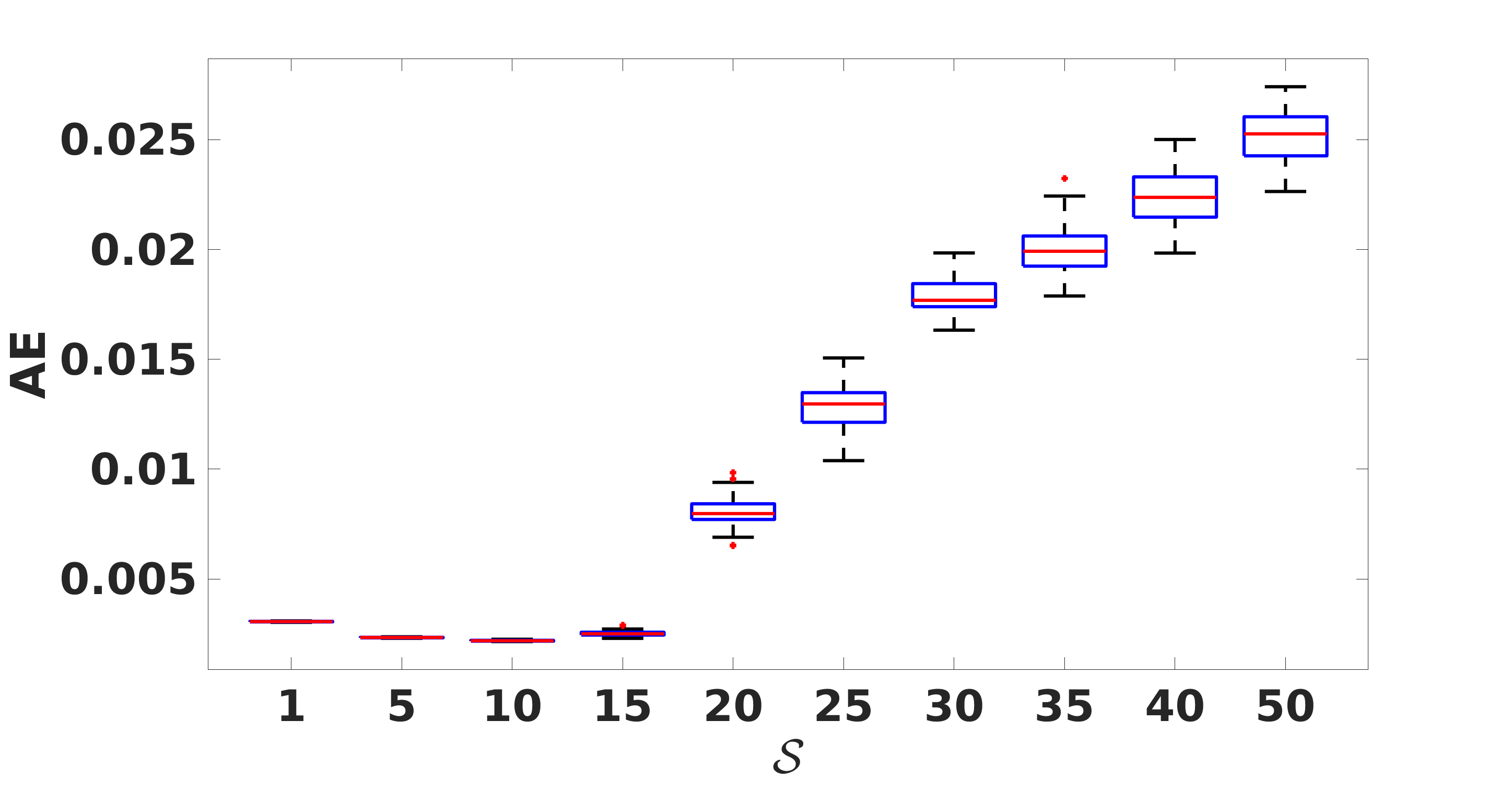}\label{part:FIG5a}} \,
    \subfloat[RRMSE]{\includegraphics[width=74mm]{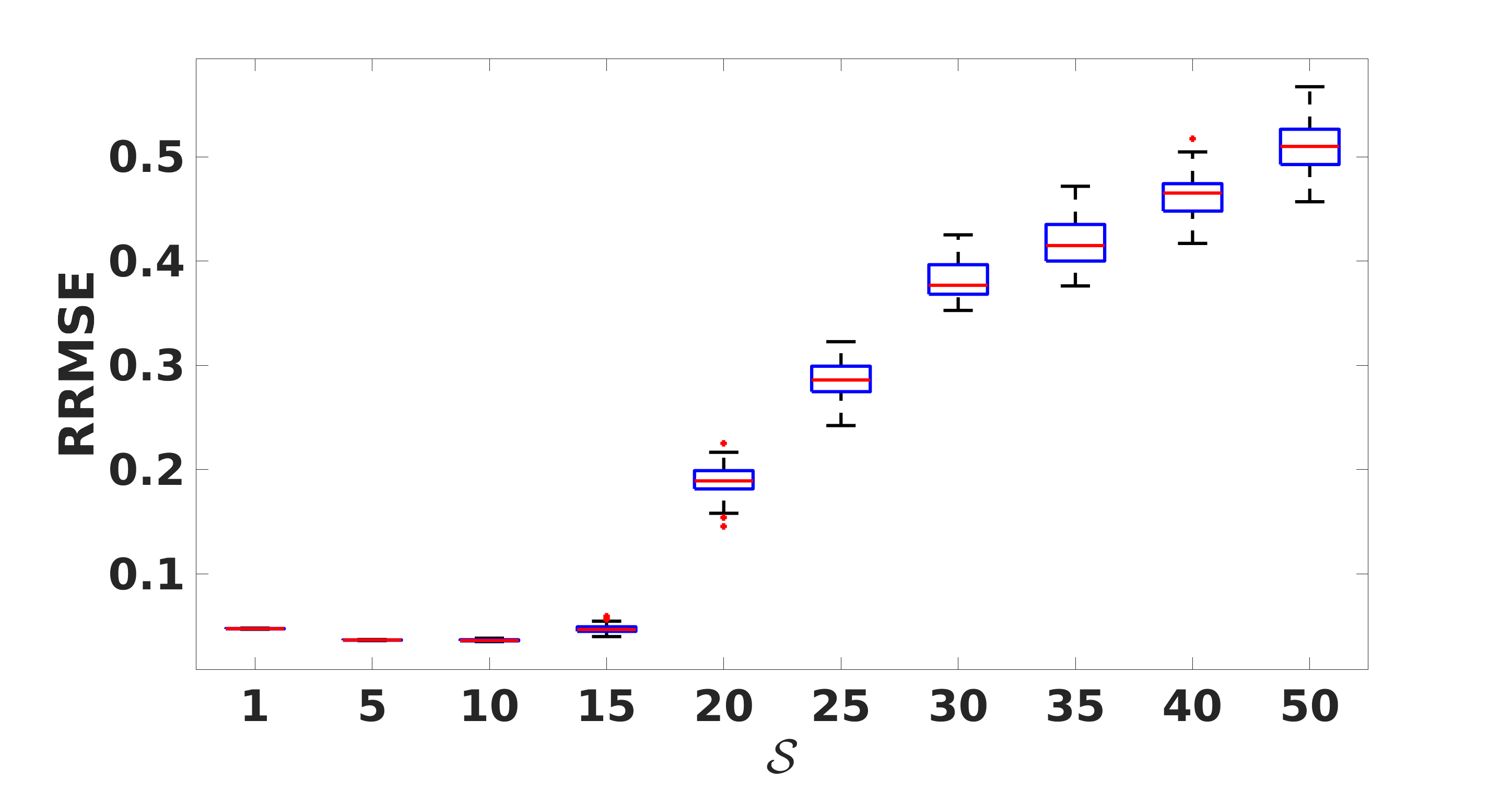}\label{part:FIG5b}} \,
    \subfloat[AES]{\includegraphics[width=74mm]{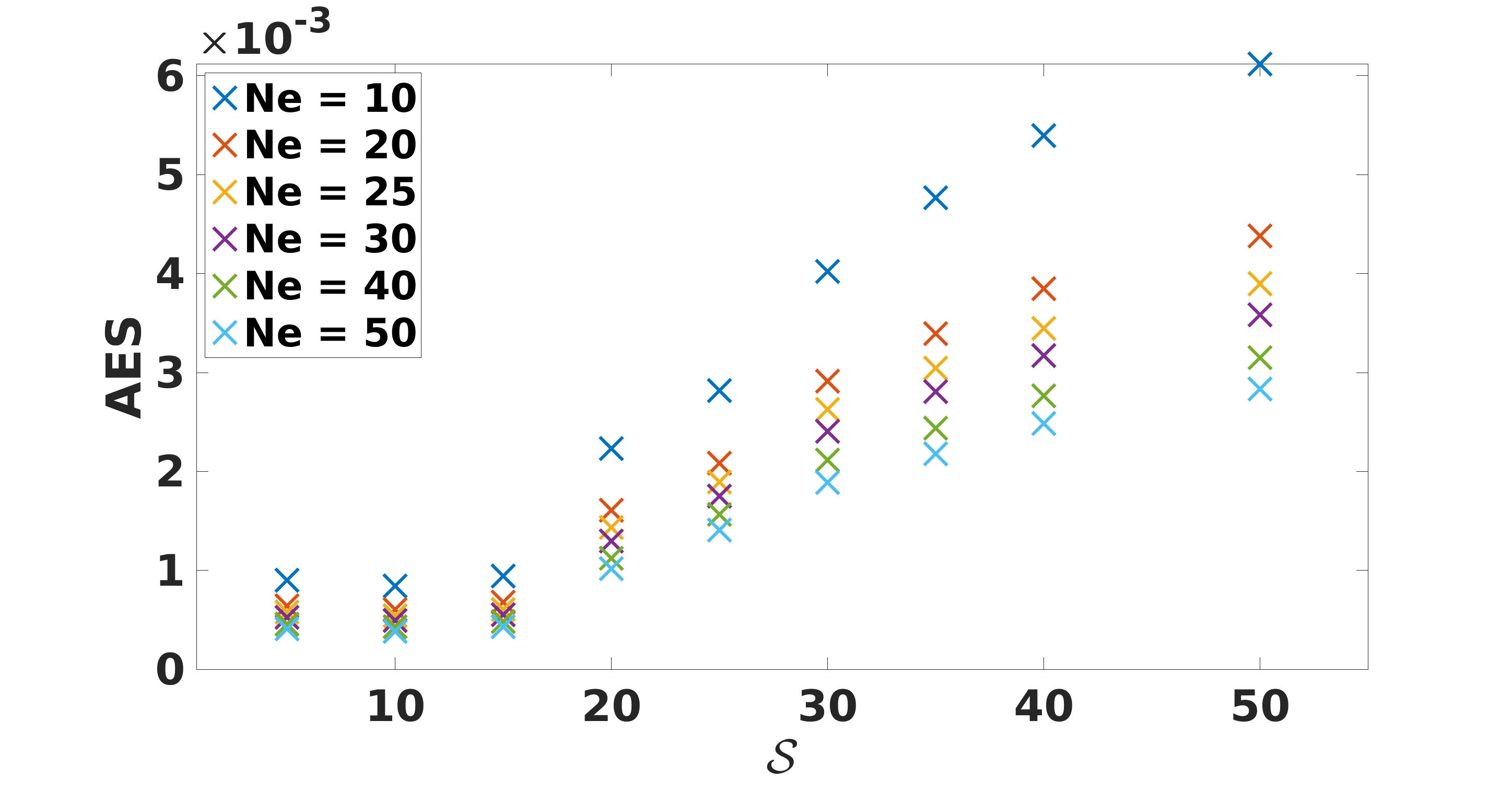}\label{part:FIG5c}} \,
    \subfloat[$\Lambda$]{\includegraphics[width=74mm]{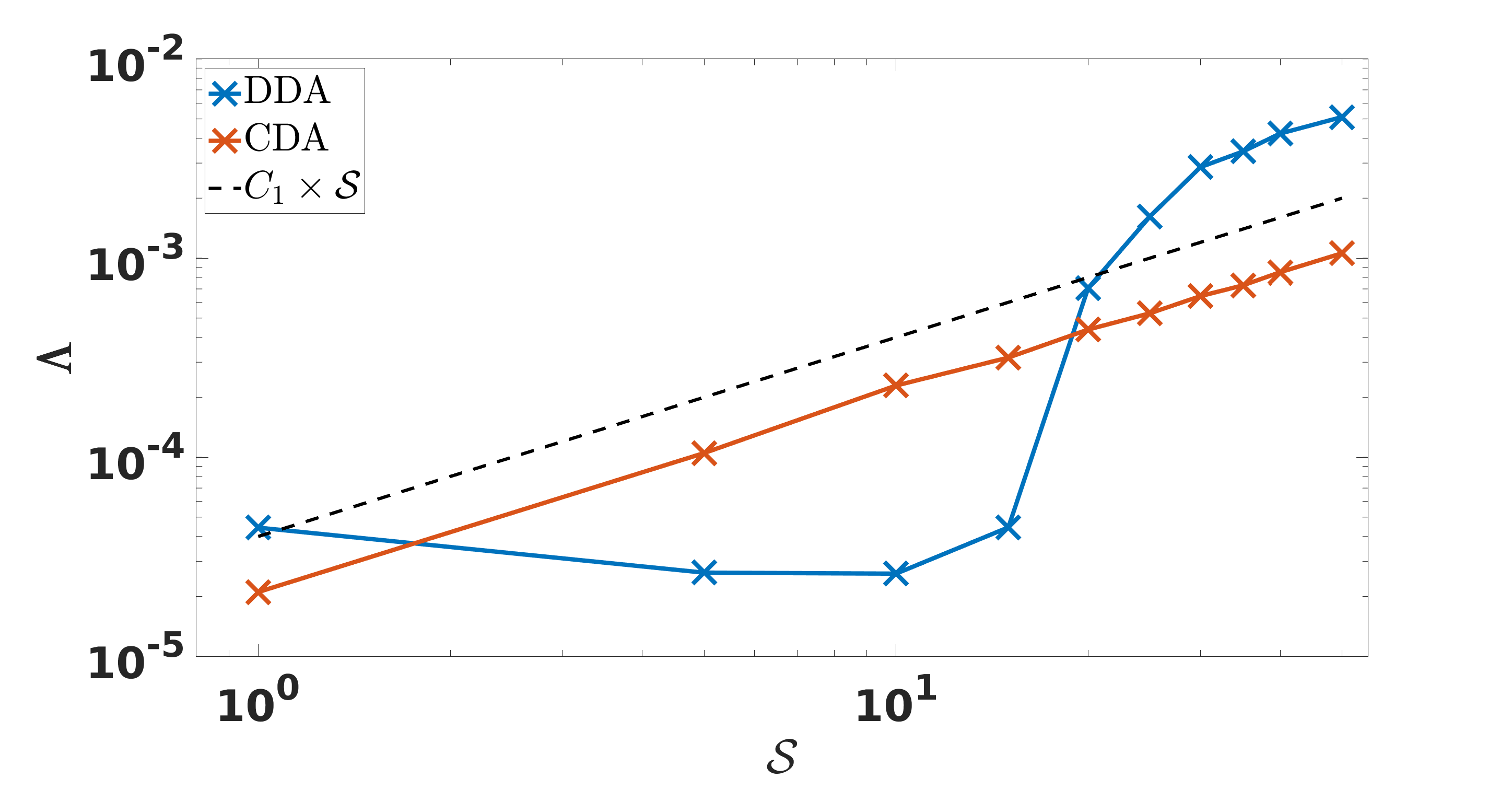}\label{part:FIG5d}} \,

\caption{Frames (a) and (b) show the box plots of the temperature skill scores, AE and RRMSE respectively, for $50$ downscaled solutions using DDA, at the final time step. Frame (c) depicts the AES for different number of realizations of the DDA downscaled fields. Frame (d) illustrates the dependence of $\Lambda$ on $\mathcal{R}$ for DDA and CDA as indicated by the legend.}
\label{lastStep_freqObs}
\end{figure}

Figure \ref{lastStep_freqObs} shows the distribution of the asymptotic values of AE, RRMSE and AES for the DDA downscaled temperature fields with varying $\mathcal{S}$ at the final time step. 
The behavior of $\Lambda$ in response to varying $\mathcal{S}$ is presented for both CDA and DDA for comparison.
Specifically, Figures \ref{part:FIG5a} and \ref{part:FIG5b} show the box plots describing the AE and RRMSE of the temperature corresponding to each of the considered observation frequencies. 
For $\mathcal{S} \leq 15$, the distributions of the AE and RRMSE are tight with little variability.
The variability of these skill scores increases as $\mathcal{S}$ increases. 
The box plots also show that both the AE and RRMSE initially decrease with increasing $\mathcal{S}$ until reaching their minimum at $\mathcal{S}=10$, after which, the errors increase nonlinearly with increasing $\mathcal{S}$. 
In particular, for $\mathcal{S} \geq 25$, both AE and RRMSE show large values indicating a substantial deviation of the downscaled solution from the reference solution. 
Figure \ref{part:FIG5c} plots the variation of AES for different $\mathcal{S}$ and number of downscaled fields. 
As expected, the results suggest that for a fixed $\mathcal{S}$, the AES decreases as $1/\sqrt{N_{ens}}$.
For a fixed given number of downscaled fields, the AES decreases to a minimum as $\mathcal{S}$ increases from $1$ up to $10$, beyond which the AES continues to increase following a similar behavior to that of the AE and RRMSE.
Finally, figure \ref{part:FIG5d} illustrates the variation of $\Lambda$ with $\mathcal{S}$ for both the CDA and DDA downscaled temperature fields. 
$\Lambda$ varies linearly with $\mathcal{S}$ for CDA and exhibits larger values than those obtained with DDA, except for extreme values of $\mathcal{S}$.
On the other hand, $\Lambda$ evolves non-monotonically with $\mathcal{S}$ for the DDA downscaled fields, and decreases to a minimum as $\mathcal{S}$ increases from $1$ up to $10$ before nonlinearly increasing with increasing $\mathcal{S}$.

The behavior of DDA errors is attributed to the discrete nature of the nudging, which allows the solution to evolve based on the exact system of equations between the nudging time steps. 
At very small values of $\mathcal{S}$, however, nudging is frequently applied, and thus, the solution becomes increasingly affected by the observation noise leading to an additional error.

\subsection{Sensitivity to the Observation's Spatial Resolution}
\label{sec:densObs}

The CDA was theoretically proven to converge if the observations are not too far apart \citep{Azouani2013}. 
Specifically, the resolution of the observational grid must fall within the range of the dissipation spectrum. 
In this experiment, the response of CDA and DDA to different spatial resolutions of the observational grid is examined. 
Downscaling runs were then performed for $\mathcal{R}$ of $5$, $15$ and $25$, all with $\mathcal{S} = 10$. 
The joint effect of coarse observations in space and time is also examined with $\mathcal{S} = 25$ and $\mathcal{R}=5$ and $15$.
The nudging factors of DDA and CDA were set to $7.0$ and $3.0$, respectively. 
The observational noise levels are set to their nominal levels, $(\sigma_T, \sigma_{\mathbf{u}}) = (0.1,0.05)$.

\begin{figure}[!htbp]
\centering
\hspace{-1cm}

\subfloat[$\mathcal{R} = 5, \mathcal{S} = 10$]{\includegraphics[width=50mm]{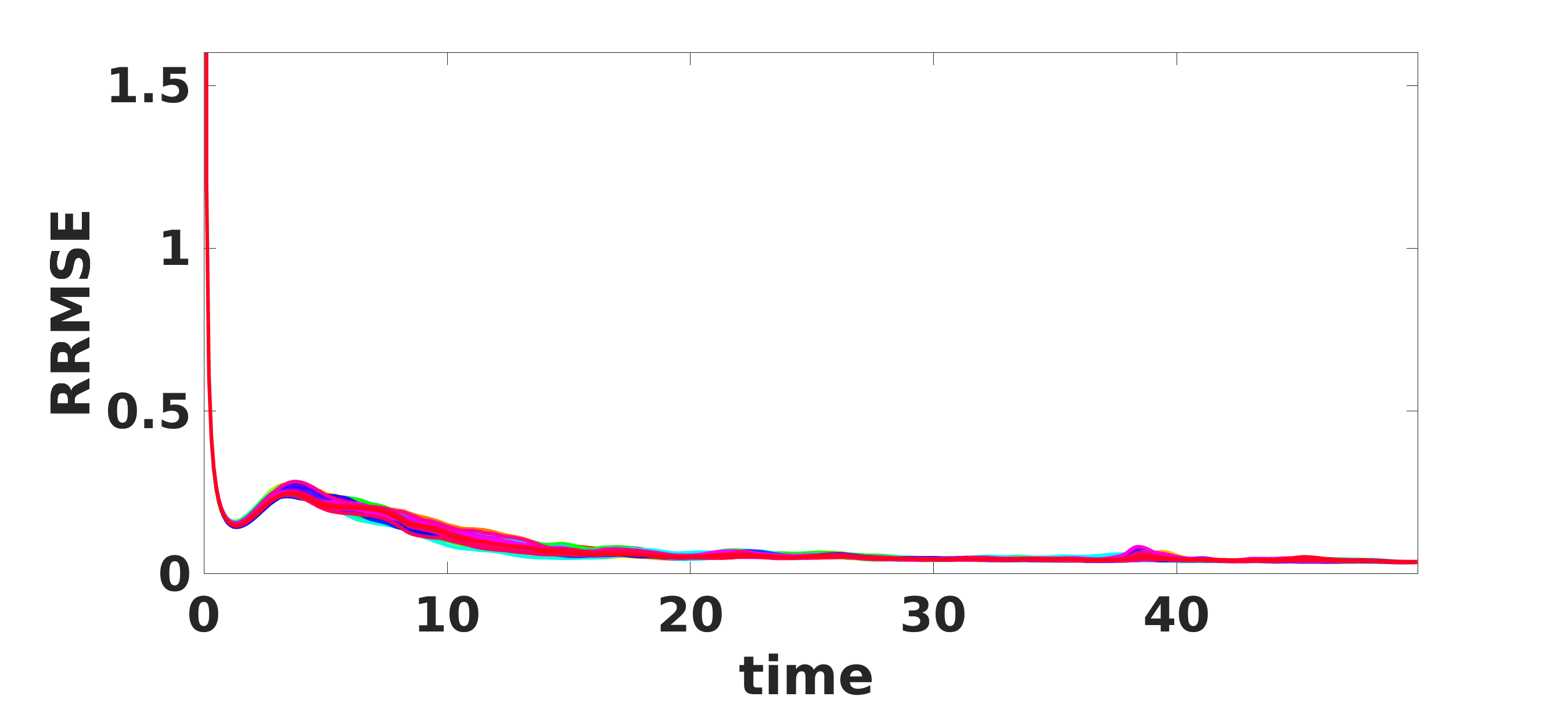}} \,
\subfloat[$\mathcal{R} = 5, \mathcal{S} = 25$]{\includegraphics[width=50mm]{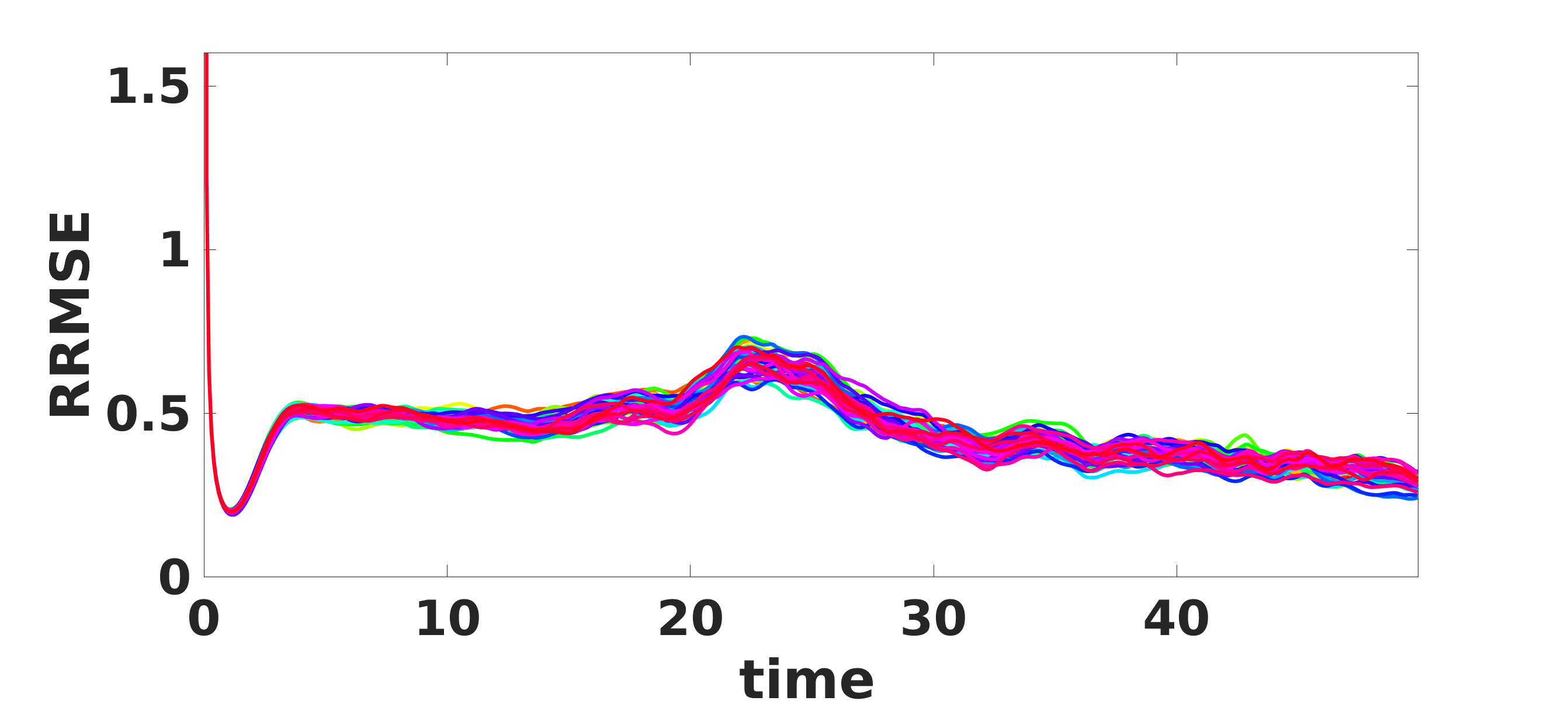}} \,
\subfloat[$\mathcal{R} = 15, \mathcal{S} = 10$]{\includegraphics[width=50mm]{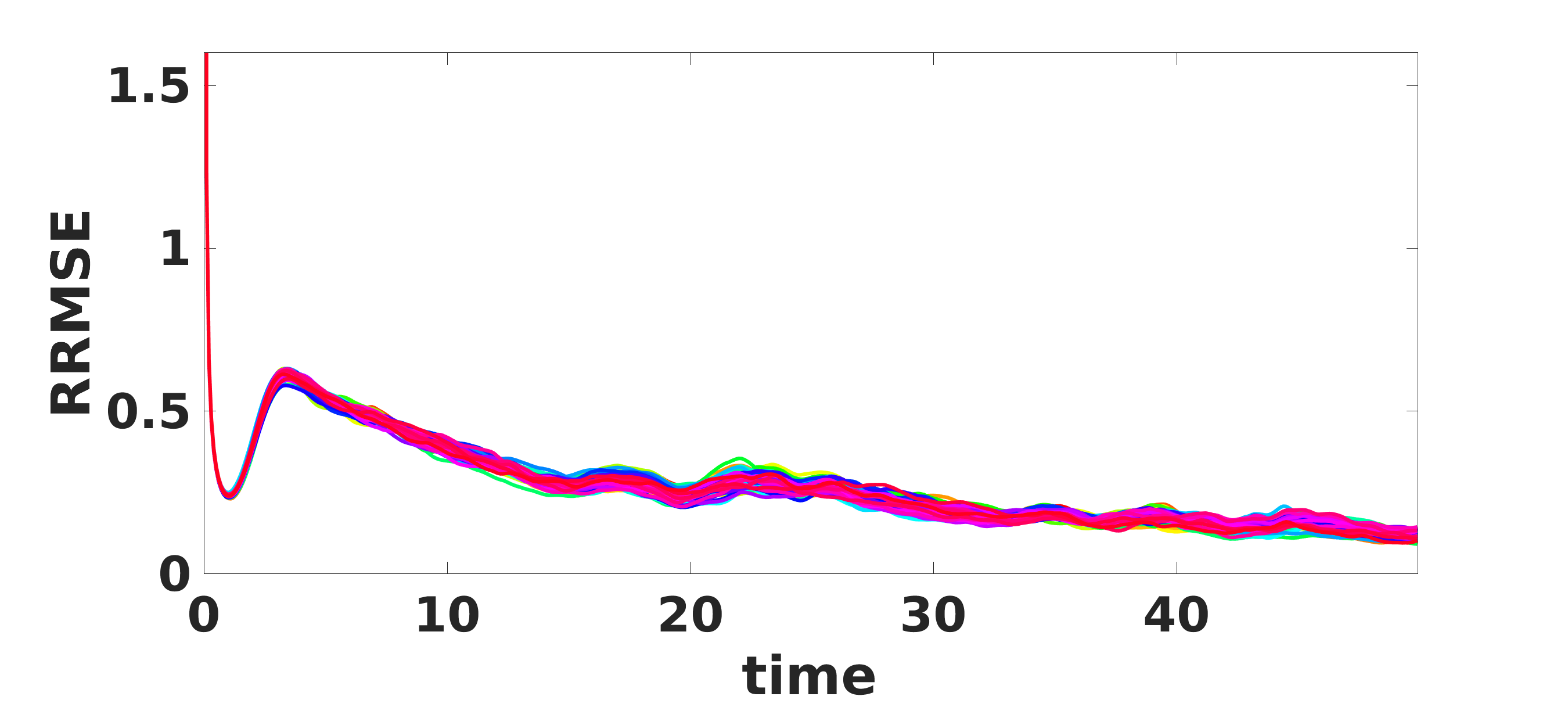}} \,
\subfloat[$\mathcal{R} = 15, \mathcal{S} = 25$]{\includegraphics[width=50mm]{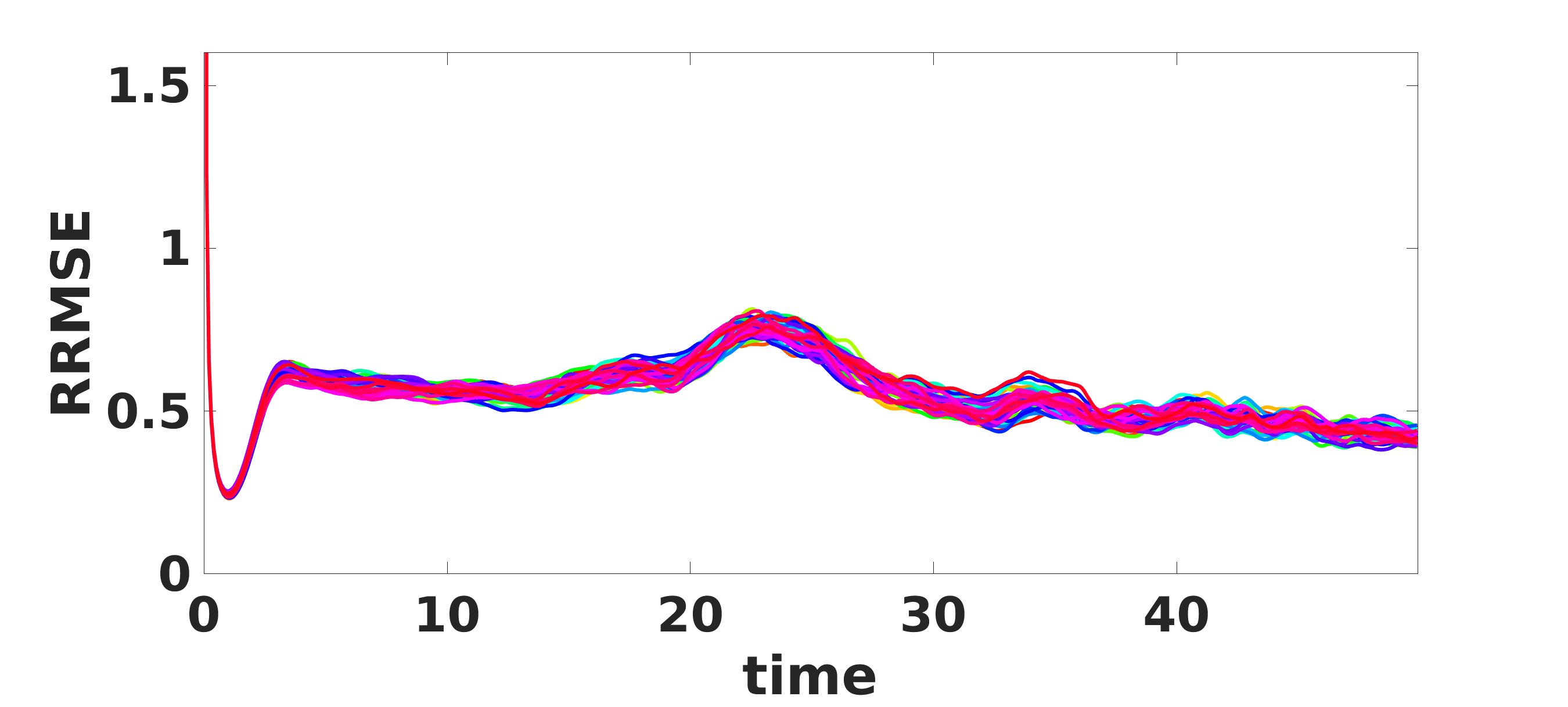}} \,
\subfloat[$\mathcal{R} = 25, \mathcal{S} = 10$]{\includegraphics[width=50mm]{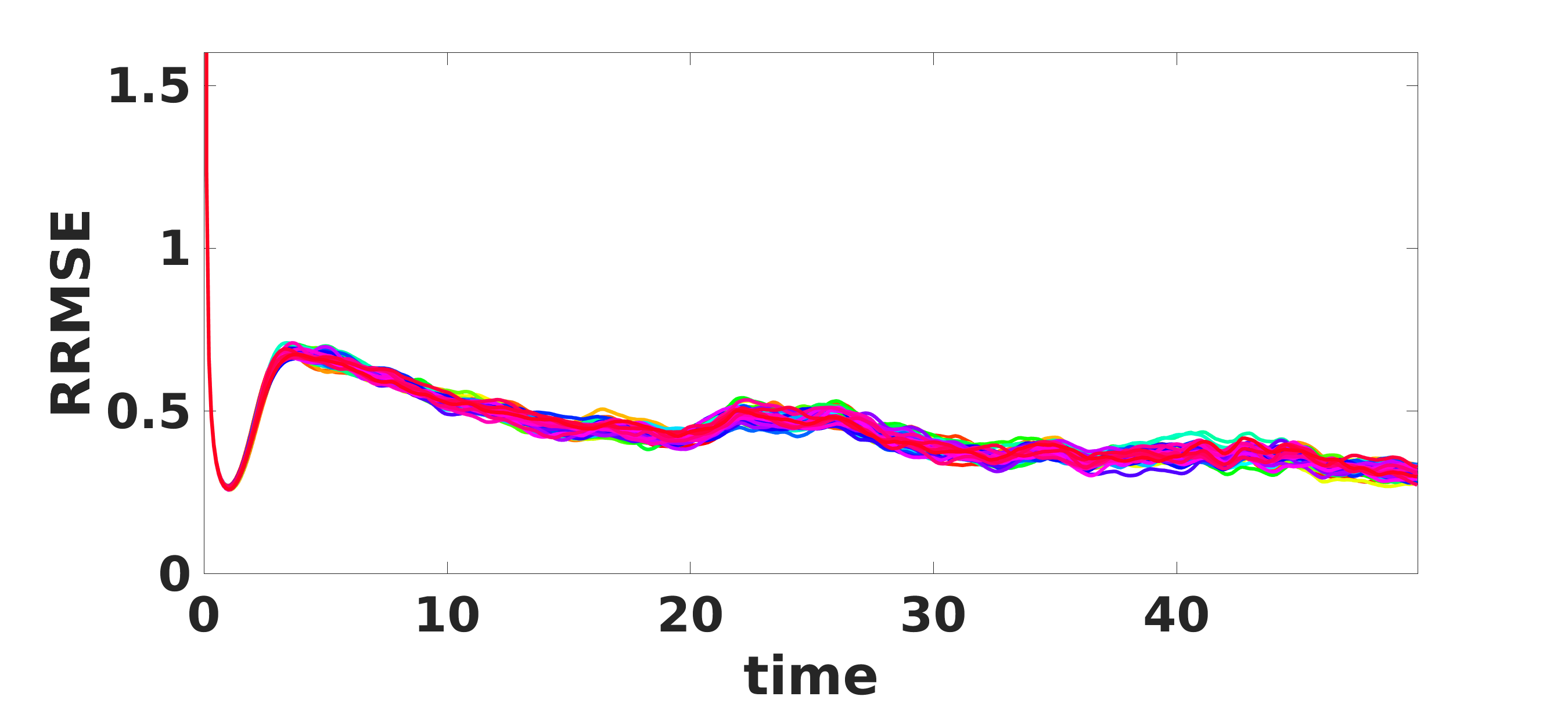}} \,

\caption{Evolution of the RRMSE of $50$ downscaled temperature fields corresponding to DDA. Plots are generated for different $\mathcal{R}$ and $\mathcal{S}$ as indicated.}
\label{evo:gridDens}
\end{figure}

Figure \ref{evo:gridDens} shows the time evolution of the temperature RRMSE for the different combinations of $\mathcal{S}$ and $\mathcal{R}$ values. 
In all cases, the RRMSE exhibits an initial rapid drop followed by an increase before slowly settling down to a plateau. 
For $\mathcal{S}=10$, the height of the plateau noticeably increases as $\mathcal{R}$ increases. 
Specifically, for $\mathcal{R}=15$, all the RRMSEs of the downscaled realizations experience an initial rapid drop followed by a small increase ending with a plateau.
As $\mathcal{R}$ increases, however, the RRMSE curves clearly exhibited a larger spread and reached a higher value at the plateau, which took a longer period to reach. 
On the other hand, for $\mathcal{S}=25$ and for both $\mathcal{R}$ values examined, the RRMSE curves exhibit a smaller rise in the plateau with increasing $\mathcal{R}$ in comparison to the case of $\mathcal{S}=10$. 
The monotonic decrease in the CDA's RRMSE is lost with increasing $\mathcal{R}$, as seen in the corresponding results in the Supplementary material. 
In particular, for $\mathcal{R}=15$ and $25$, the RRMSE rapidly decreases, then increases for a short period before decreasing to a final plateau. 
The DDA seems to be more sensitive to $\mathcal{R}$ for the smaller value of $\mathcal{S}$ indicated by the increase in the value of the plateau. 
These results suggest that the coarser the spatial resolution of the observational grid is, the less accurate the downscaled solution becomes. 
A thorough assessment of the behavior of the skill scores at the final time step is presented in the following.

\begin{figure}[!htbp]
\centering
\hspace{-1cm}

\subfloat[AE]{\includegraphics[width=74mm]{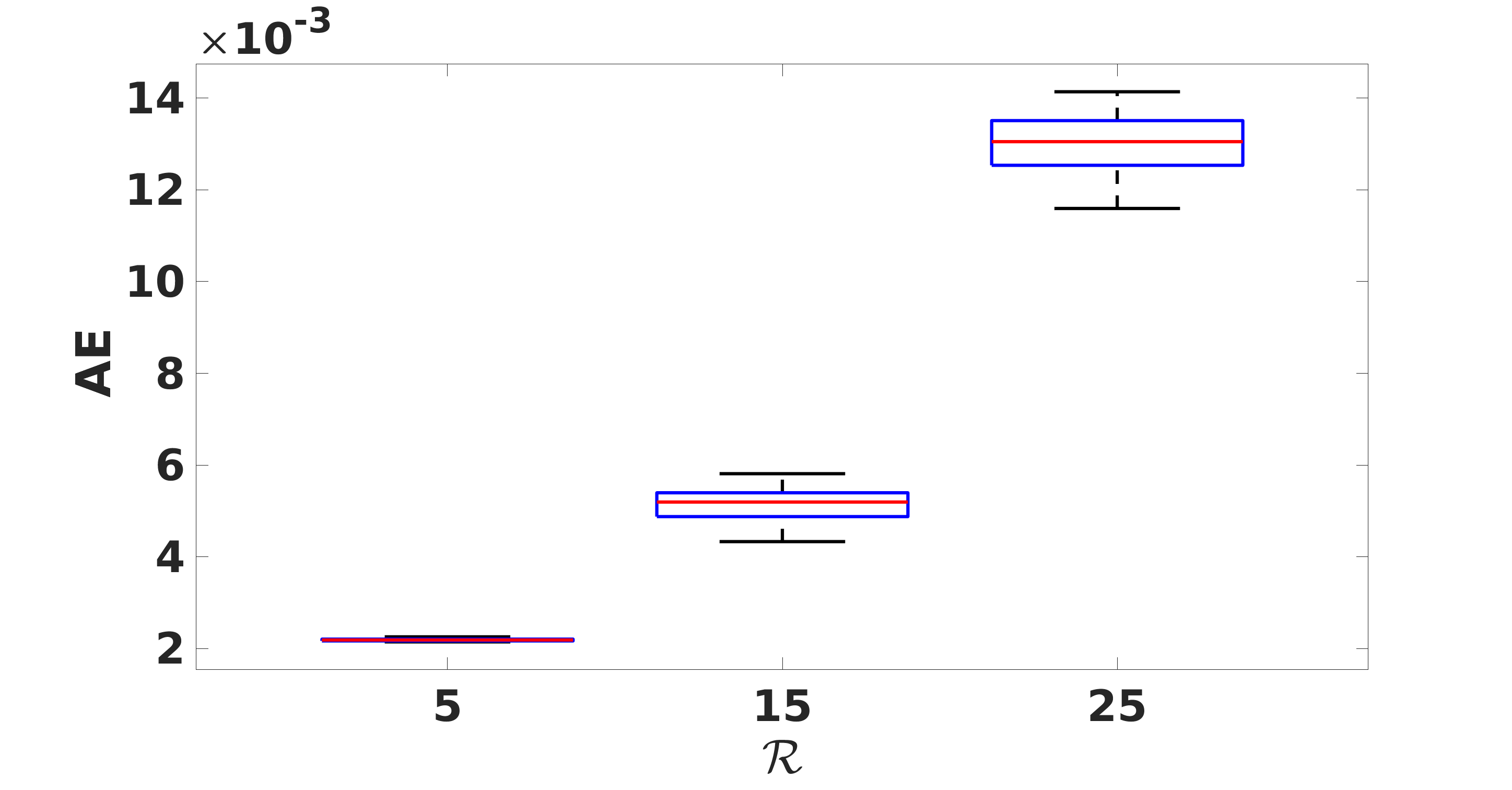}\label{7a}} \,
\subfloat[RRMSE]{\includegraphics[width=74mm]{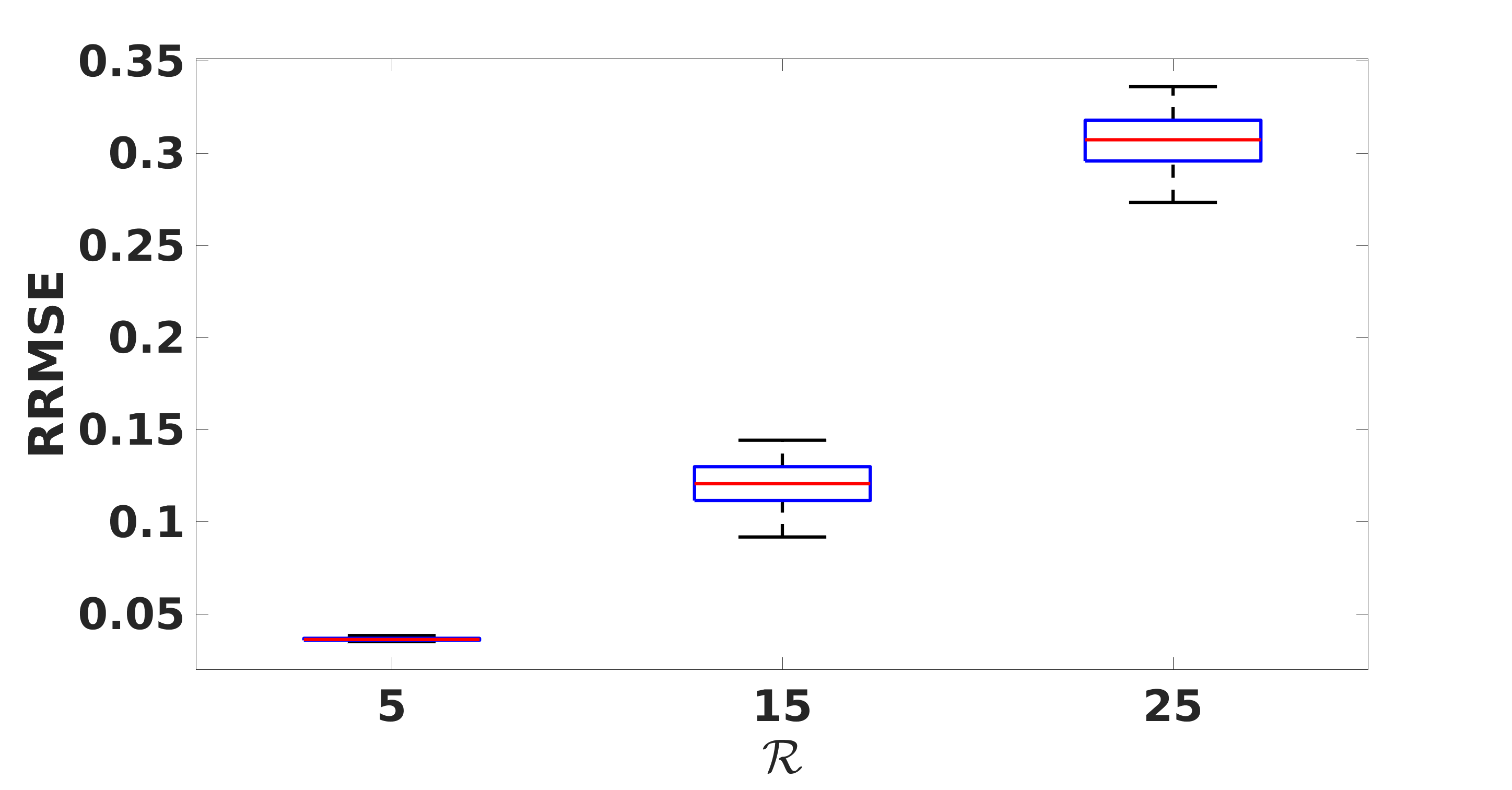}\label{7b}} \,
\subfloat[AES]{\includegraphics[width=74mm]{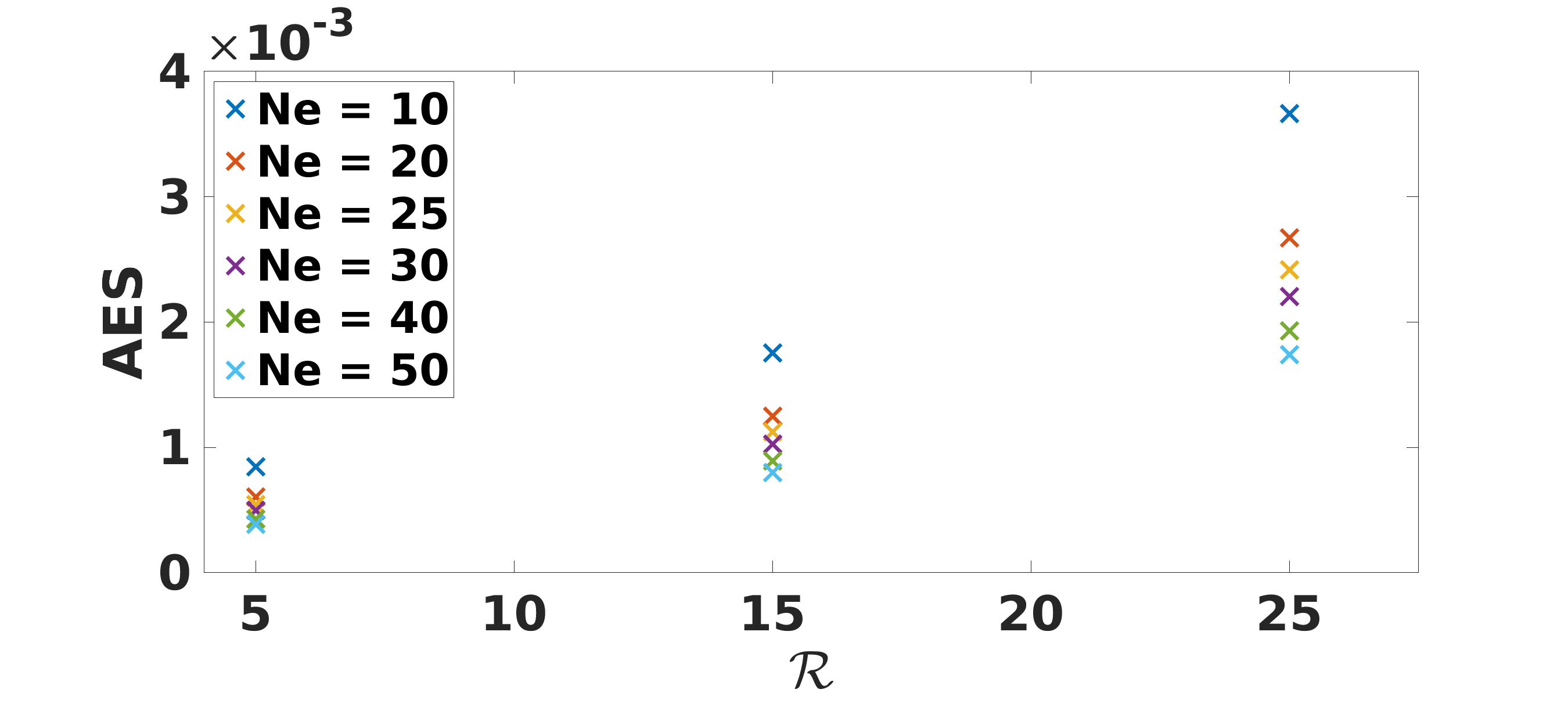}\label{7c}} \,
\subfloat[$\Lambda$]{\includegraphics[width=74mm]{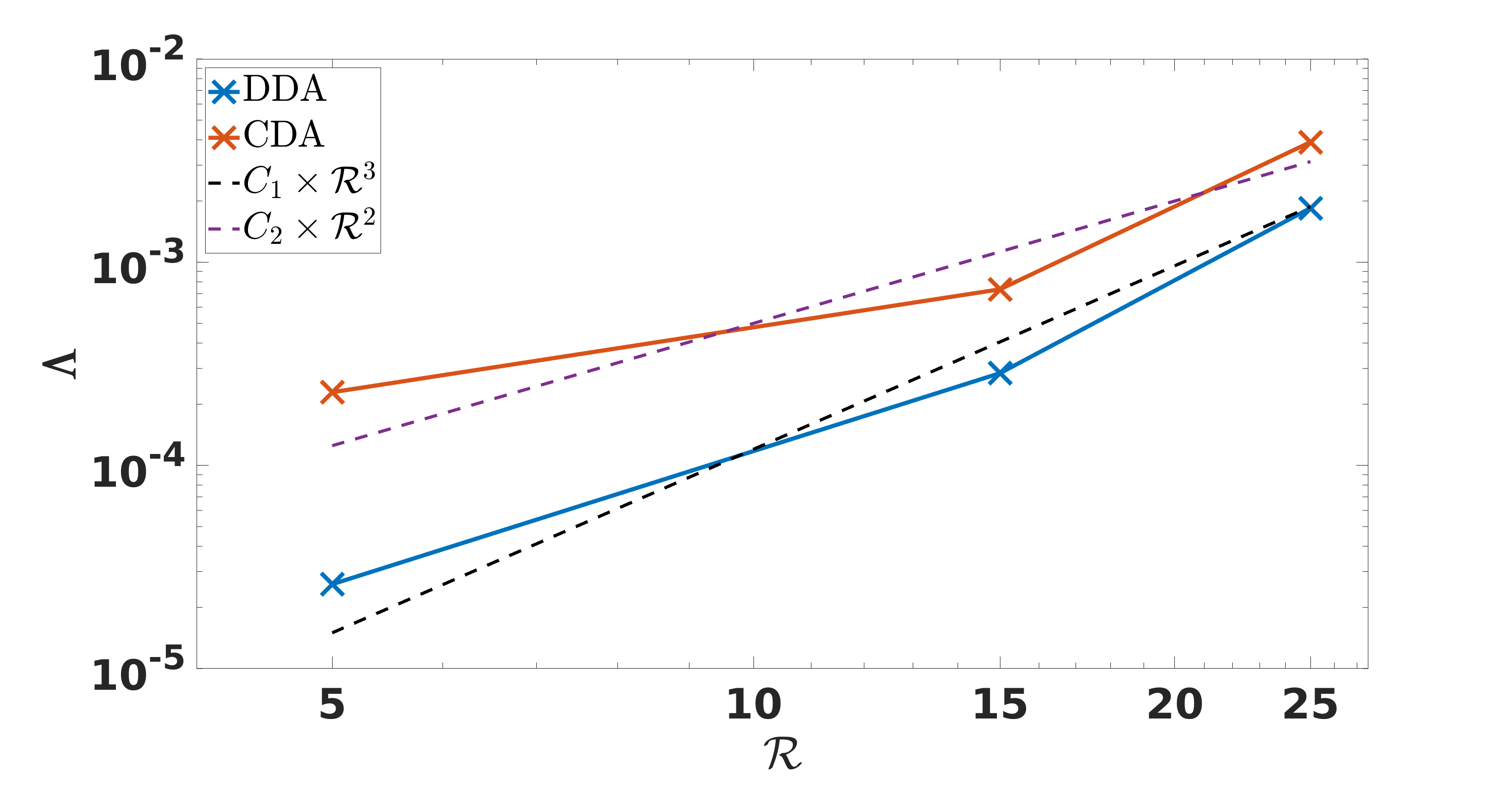}\label{7d}} \,

\caption{Frames (a) and (b) show the box plots of the temperature skill scores, AE and RRMSE respectively, for $50$ downscaled solutions using DDA with $\mathcal{S}=10$ and the indicated $\mathcal{R}$, at  the last time step. Frame (c) depicts the AES for different number of realizations. Frame (d) illustrates the dependence of $\Lambda$ on $\mathcal{R}$ for both DDA and CDA.}
\label{lastStep:gridDensev11}
\end{figure}

\begin{figure}[!htbp]
\centering
\hspace{-1cm}

\subfloat[AE]{\includegraphics[width=90mm]{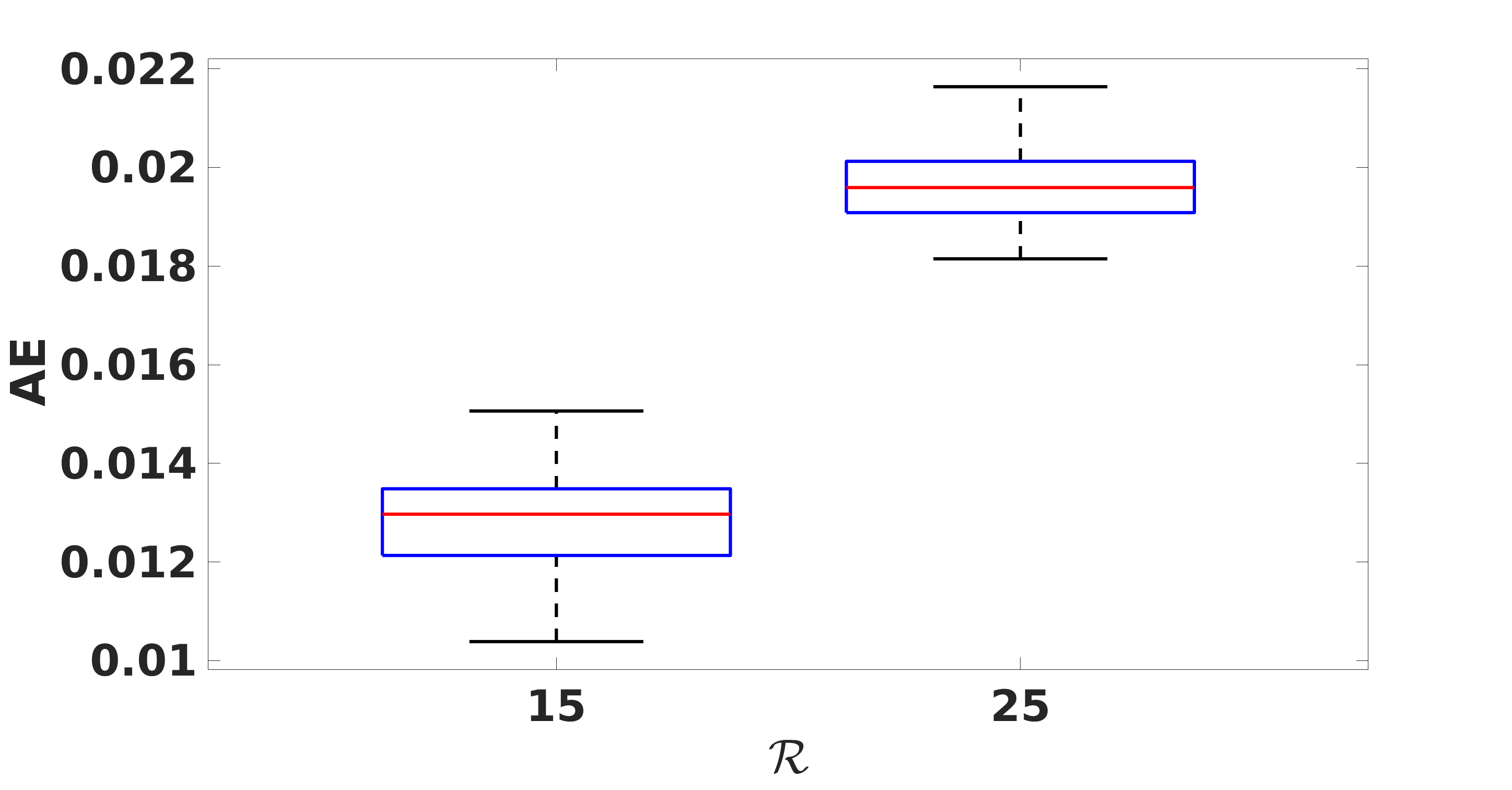} \label{8a}} \,
\subfloat[RRMSE]{\includegraphics[width=90mm]{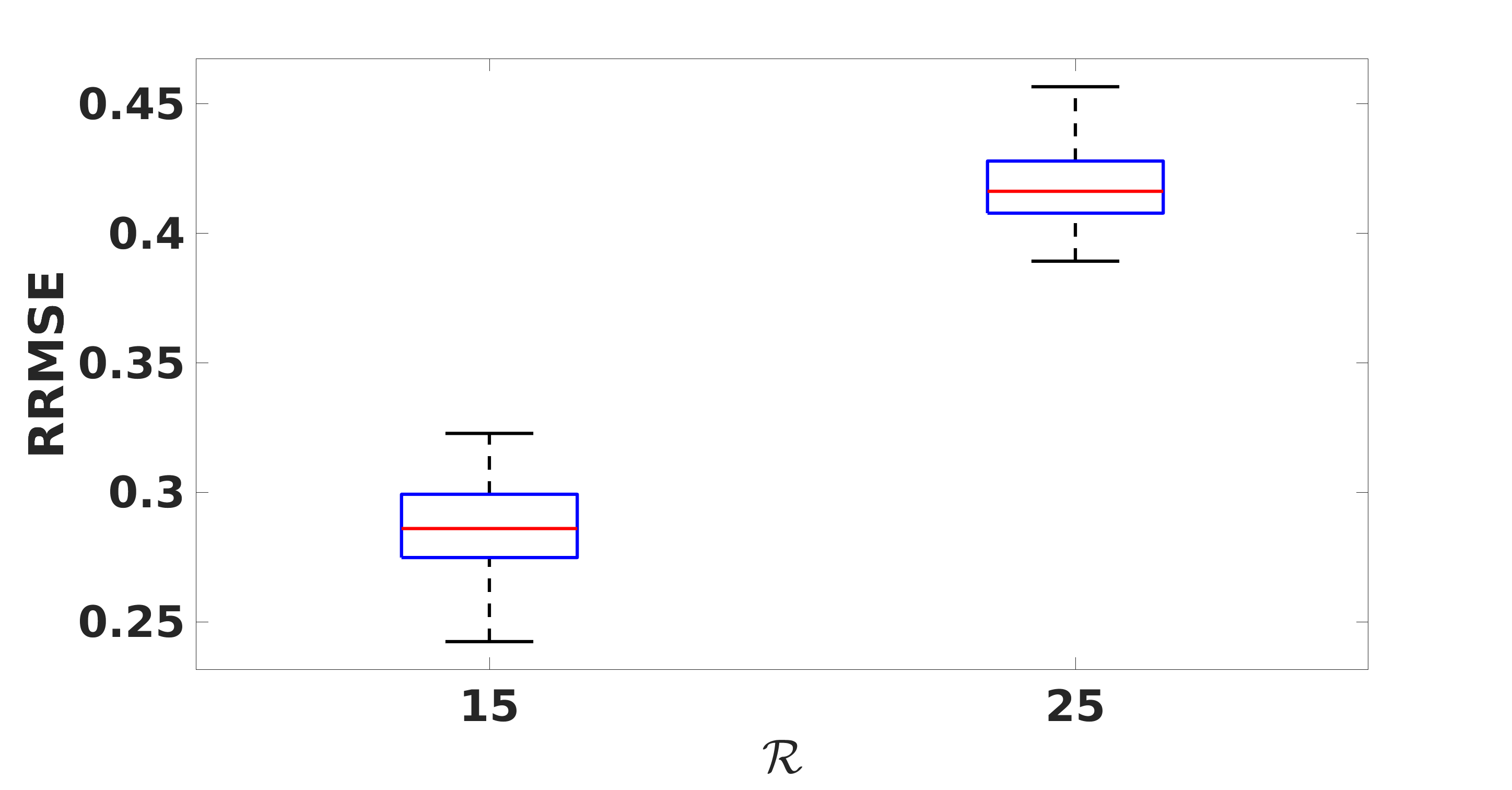} \label{8b}} \,

\caption{Box plots of the temperature skill scores: (a) AE and (b) RRMSE for $50$ downscaled solutions using DDA with $\mathcal{S}=25$ and the indicated $\mathcal{R}$}
\label{lastStep:gridDensev25}
\end{figure}

Figure \ref{lastStep:gridDensev11} presents the box plots of AE, RRMSE and AES for $\mathcal{S}=10$ and the different tested $\mathcal{R}$ values for the DDA downscaled temperature fields at the last time step. 
The figure also displays the behavior of $\Lambda$ for both CDA and DDA downscaled temperature fields at the last time step.
The plots suggest a cubic increase of the means of both the AE and RRMSE and an increase in the solutions' spread with increasing $\mathcal{R}$. 
The AES increases linearly with increasing $\mathcal{R}$ and decreases with the square-root of the number of downscaled realizations. 
Finally, Figure \ref{7d} shows the variation of $\Lambda$ as a function of $\mathcal{R}$ for both the CDA and DDA downscaled fields. 
Specifically, $\Lambda$ increases  with increasing $\mathcal{R}$ for both algorithms. 
$\Lambda$ increases quadratically with $\mathcal{R}$ with CDA, whereas with DDA, $\Lambda$ is a cubic function of $\mathcal{R}$.
For comparison, Figure \ref{lastStep:gridDensev25} shows the box plots of the AE and RRMSE at the final time step for the case with $\mathcal{S}=25$ and the considered values of $\mathcal{R}$. 
Both AE and RRMSE increase when $\mathcal{S}$ is larger and increase for increasing $\mathcal{R}$. 
In comparison to the results of Figure \ref{lastStep:gridDensev11}, the mean of the asymptotic errors almost triple, for a fixed $\mathcal{R}$, when $\mathcal{S}$ is larger. 
In addition, the variance of the downscaled fields becomes larger with larger $\mathcal{S}$.

\subsection{Distribution of the Downscaled Fields}
\label{sec:probDist}

Here, we analyze the distribution of the downscaled fields by applying a statistical test. 
In this experiment, the noise level was set to the nominal value, the DDA and CDA nudging parameters were fixed to $7.0$ and $3.0$, respectively, with $\mathcal{S}$ and $\mathcal{R}$ equal to $11$ and $5$, respectively. 
$200$ solution samples were then generated and their profiles in the planes $y=0.25$, $0.50$ and $0.75$ were selected. 
The Kolmogorov-Smirnov test was then applied to the solution points in these planes at discrete times to assess their distribution. 
The goal is to test for Gaussianity of the downscaled solutions in space at different instances in time.

\begin{figure}[!htbp]
\centering
\hspace{-1cm}

\subfloat[Hypothesis Test, t = $30$]{\includegraphics[width=74mm]{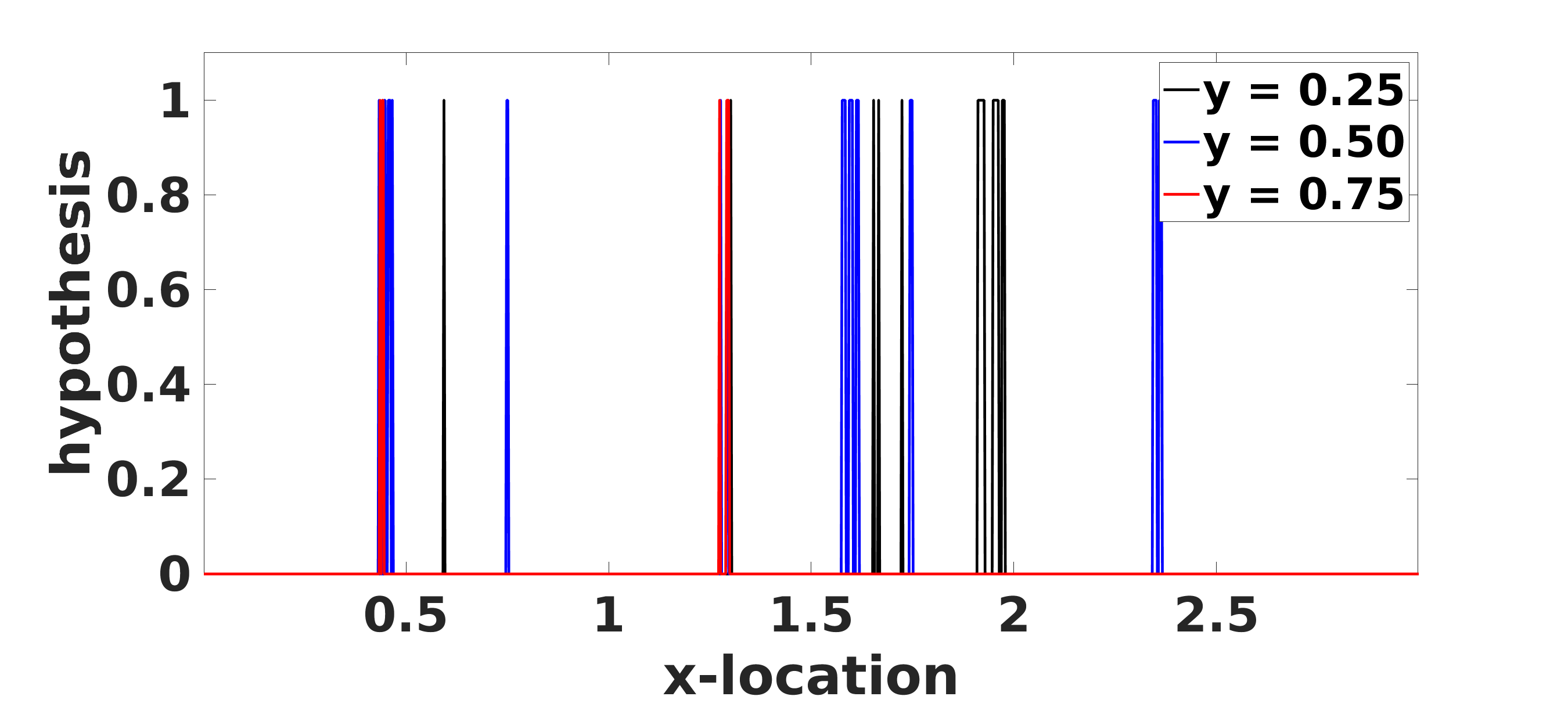}} \,
\subfloat[Hypothesis Test, t = $49.9$]{\includegraphics[width=74mm]{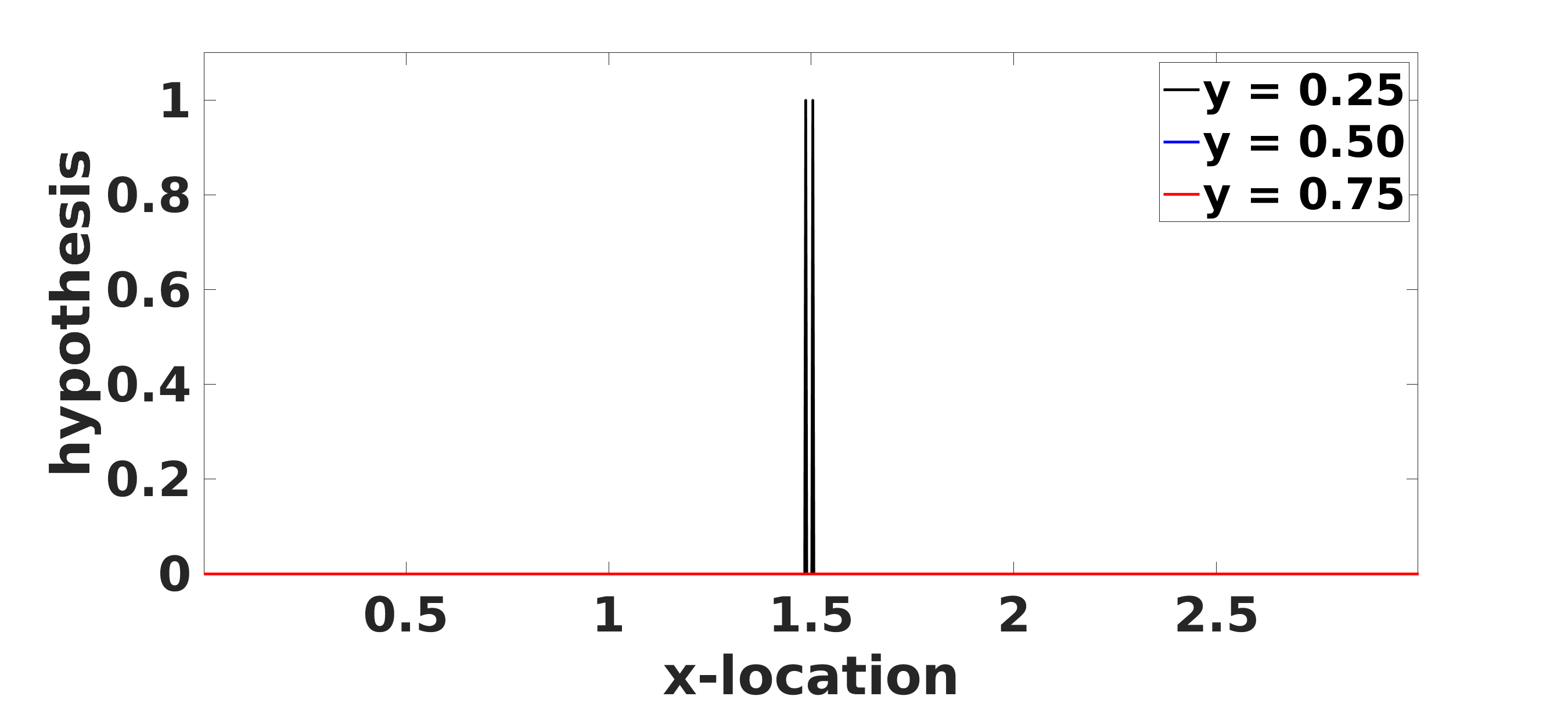}} \,
\subfloat[P-Value, t = $30$]{\includegraphics[width=74mm]{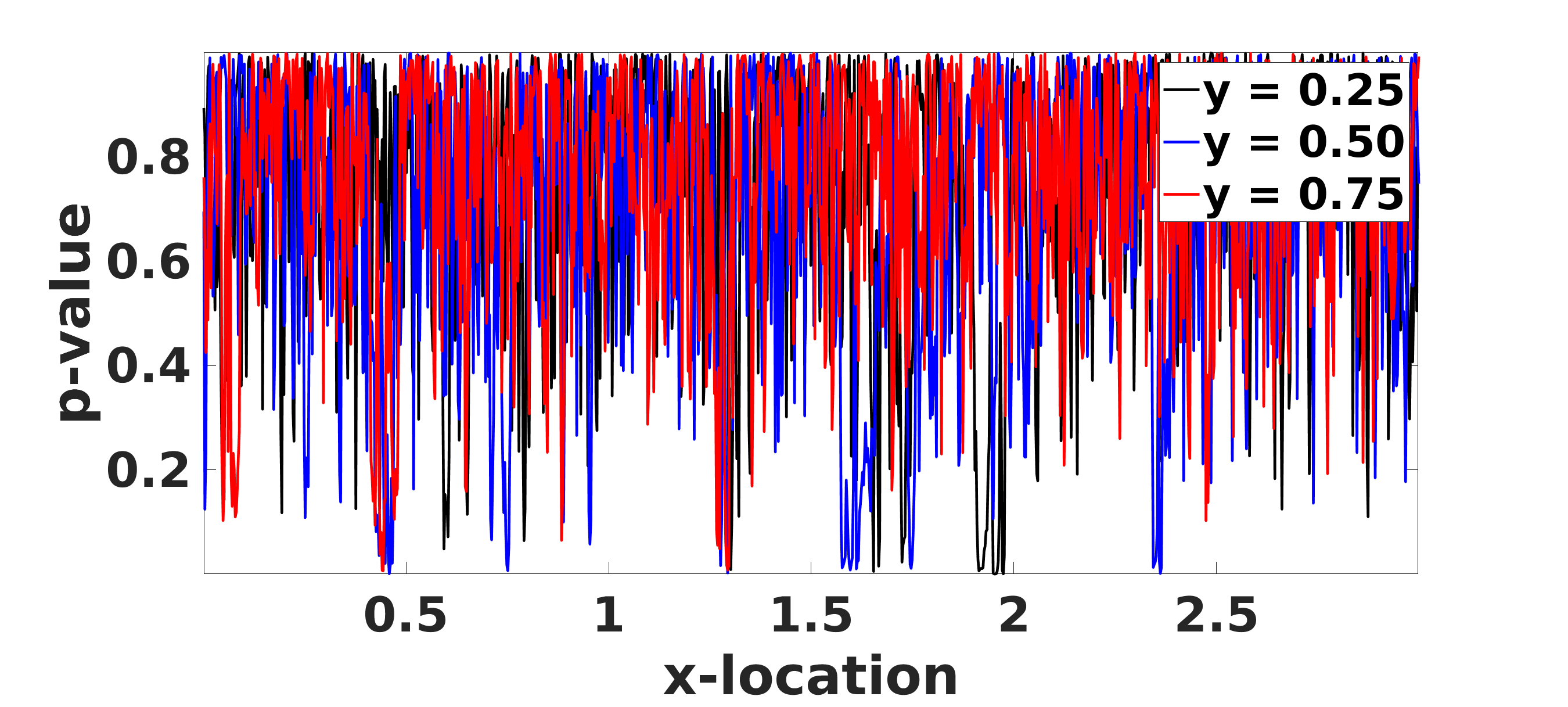}} \,
\subfloat[P-Value, t = $49.9$]{\includegraphics[width=74mm]{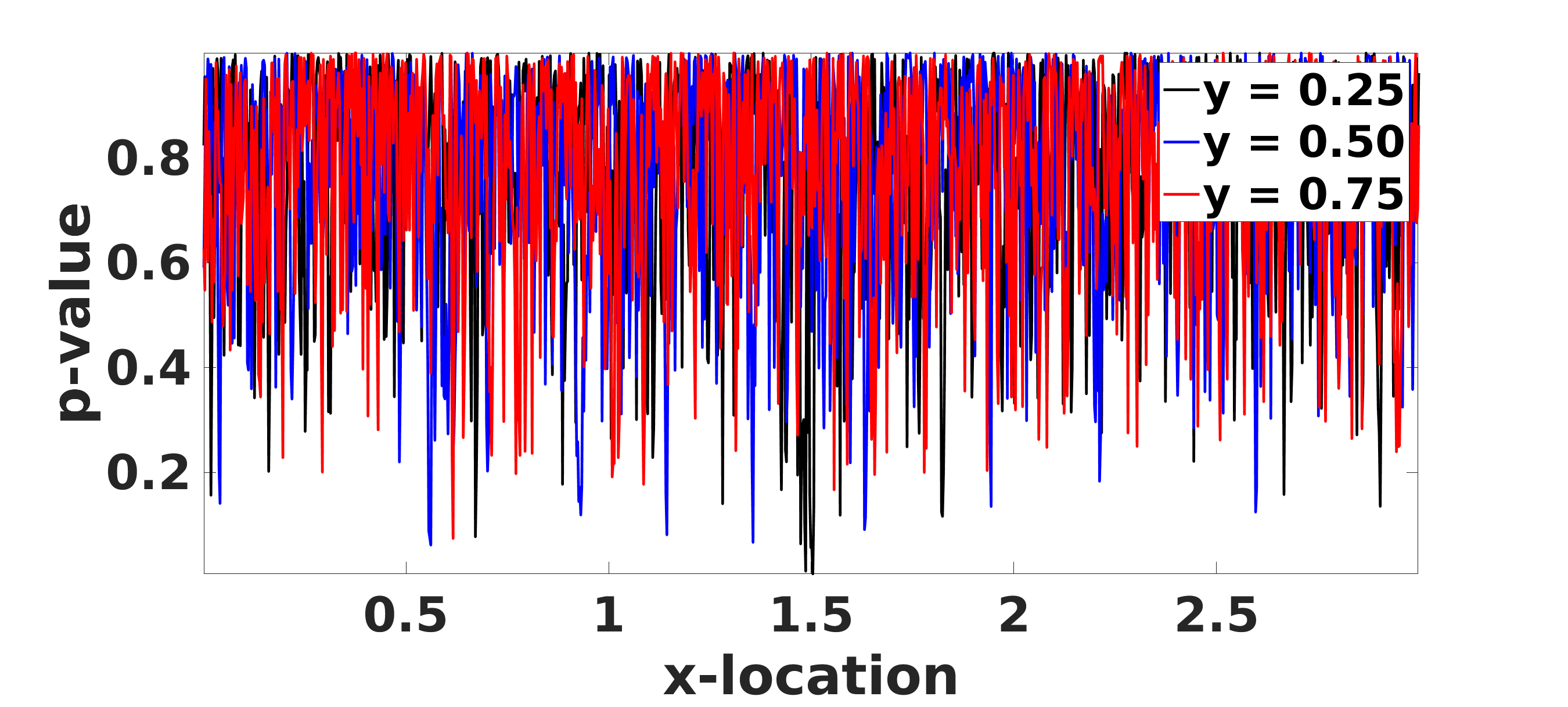}} \,

\caption{Plots of the spatial distribution of the hypothesis value and p-value for the Kolmogorov-Smirnov test. In each plot, distributions are generated based on horizontal profiles at $y=0.25$, $y=0.5$ and $y=0.75$ as indicated. Plots are generated at $t=30$ and $t=49.9$ using a $200$ realizations of the downscaled temperature solutions.}
\label{fig:gaussianity}
\end{figure}

Figure \ref{fig:gaussianity} plots snapshots of the p-value and hypothesis test for $t=30$ and $t=49.9$, representing an intermediate and the final time steps. 
The plots for $t=30$ show several locations where the hypothesis tests positive, hypothesis is equal to $1$ or equivalently p-value less than $0.05$, meaning that the samples fail the Gaussianity test. 
However, the plots at $t=49.9$ indicate that only two out of the total $3600$ tested locations fail the Gaussianity test, suggesting that, for iid Gaussian noise, the downscaled solutions are normally distributed once the plateau is reached. 
Similar results were obtained for the CDA downscaled samples, but those passed the Gaussianity test for all locations at an earlier time due to its faster convergence rate.

\subsection{Precision of the Skill Scores and a Lower Error Solution}
\label{sec:bootstrap}

In these experiments, we first assess the precision of the estimate of the RRMSE by means of a bootstrap methodology. 
Second, we explore the potential of using the average of the downscaled solutions as the best solution.
To assess the precision of the skill score estimates, bootstrap sampling is applied to subsets of various sizes of the ensemble of $200$ downscaled solutions of the previous section.
A total of $500$ bootstrap subsets were generated for each skill score.

The average of each bootstrapped skill score subset is considered as a bootstrap estimate of the skill score.
These samples were then input to a Kernel Density Estimator (KDE) to generate the probability density function (PDF) of the skill scores.

\begin{figure}[!htbp]
\centering
\hspace{-1cm}

\subfloat[Distribution - RRMSE]{\includegraphics[width=100mm]{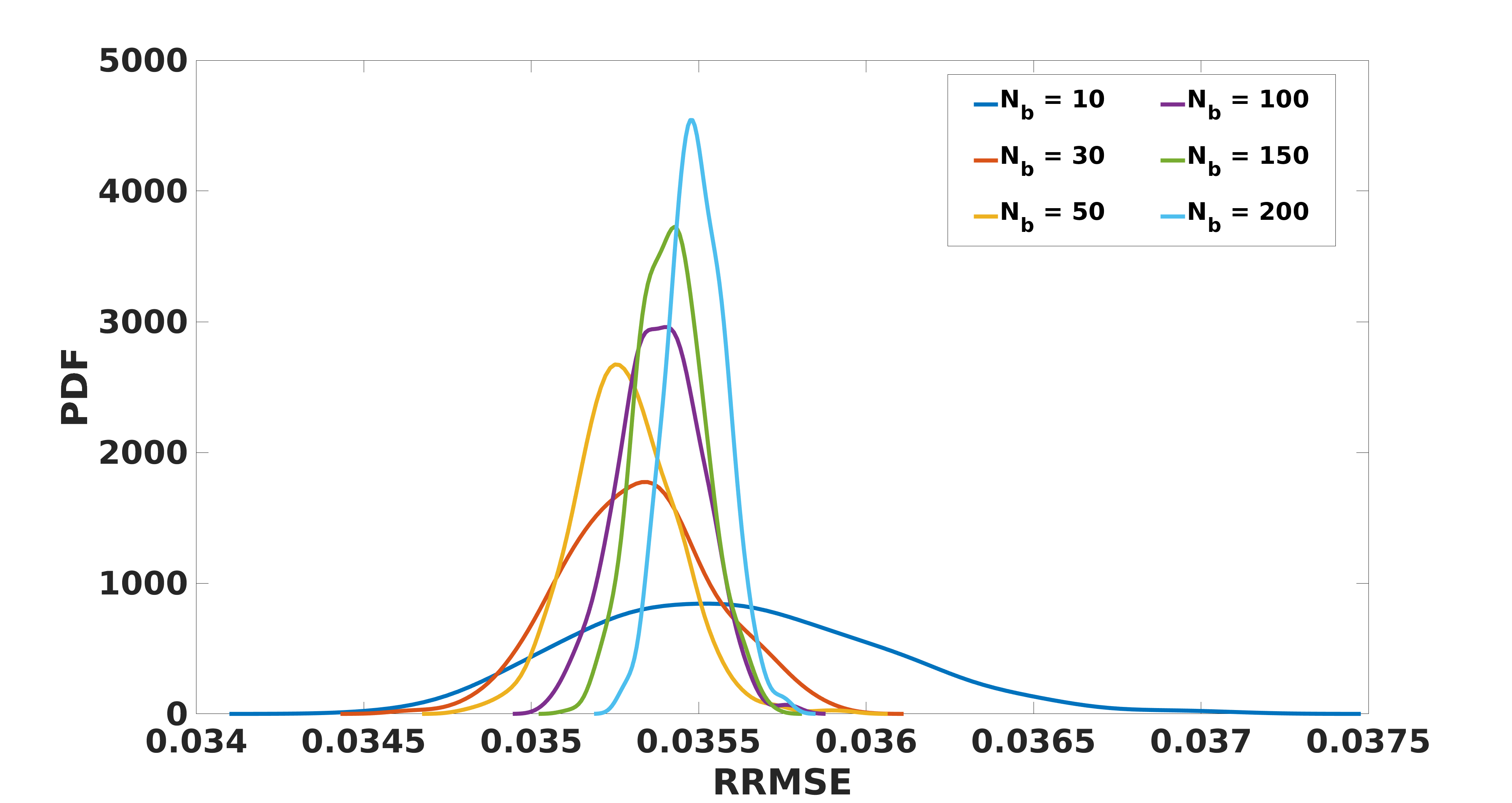}\label{10a}} \,
\subfloat[Boxplot - RRMSE]{\includegraphics[width=100mm]{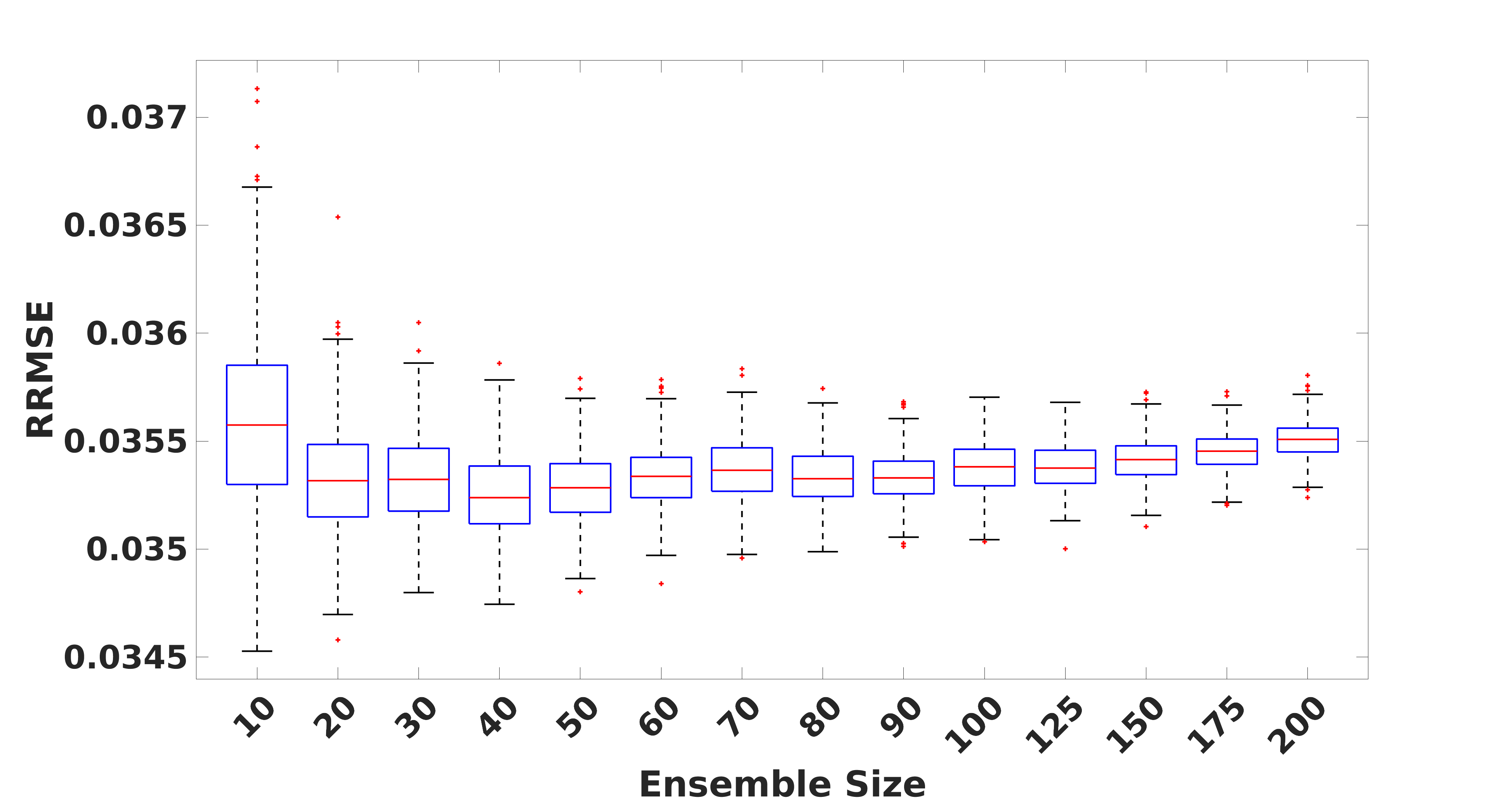}\label{10b}} \,

\caption{Distribution of the RRMSE of the bootstrapped temperature solutions for different number of realizations, as indicated, at the last time step. Plot (a) shows the PDF of the bootstrapped RRMSE of the temperature, generated using KDE, for different number of realizations as indicated, and (b) the box plots representing the the distribution of the $500$ boostrap RRMSE samples for different number of realizations as indicated.}
\label{fig:bootstrapEst}
\end{figure}

Figure \ref{fig:bootstrapEst} presents the resulting distributions for subsets of different sizes. 
Specifically, Figure \ref{10a} illustrates the distribution of the average RRMSE estimate for the different subsets and suggests that the variance in the estimated average RRMSE decreases with subset size.
The box plots describing the distribution of the boostrapped RRMSE are illustrated in Figure \ref{10b}. 
The figure shows a shrinking variability of the estimated average RRMSE with increasing subset size. 
In particular, by bootstrapping a subset with as few as $20$ realizations, the standard deviation of the distribution of the skill scores estimate can be approximately halved in comparison to the subset of $10$ realizations. 
Moreover, beyond $30$ solution samples, the standard deviation of the average RRMSE asymptotes to a steady value.

\begin{figure}[!htbp]
\centering
\hspace{-1cm}

\subfloat[Boxplot - RRMSE]{\includegraphics[width=130mm]{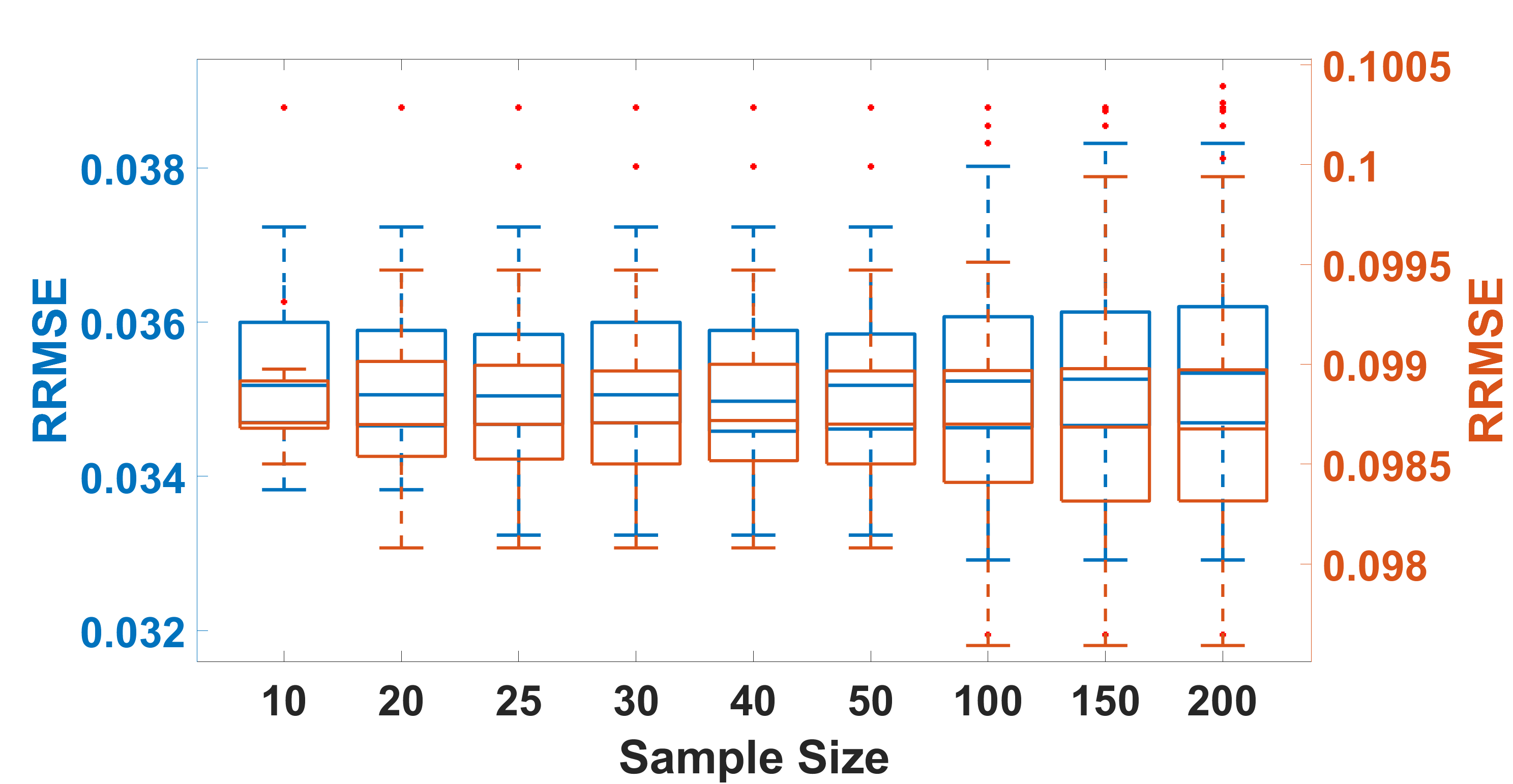}} \,
\subfloat[RRMSE vs Sample size]{\includegraphics[width=130mm]{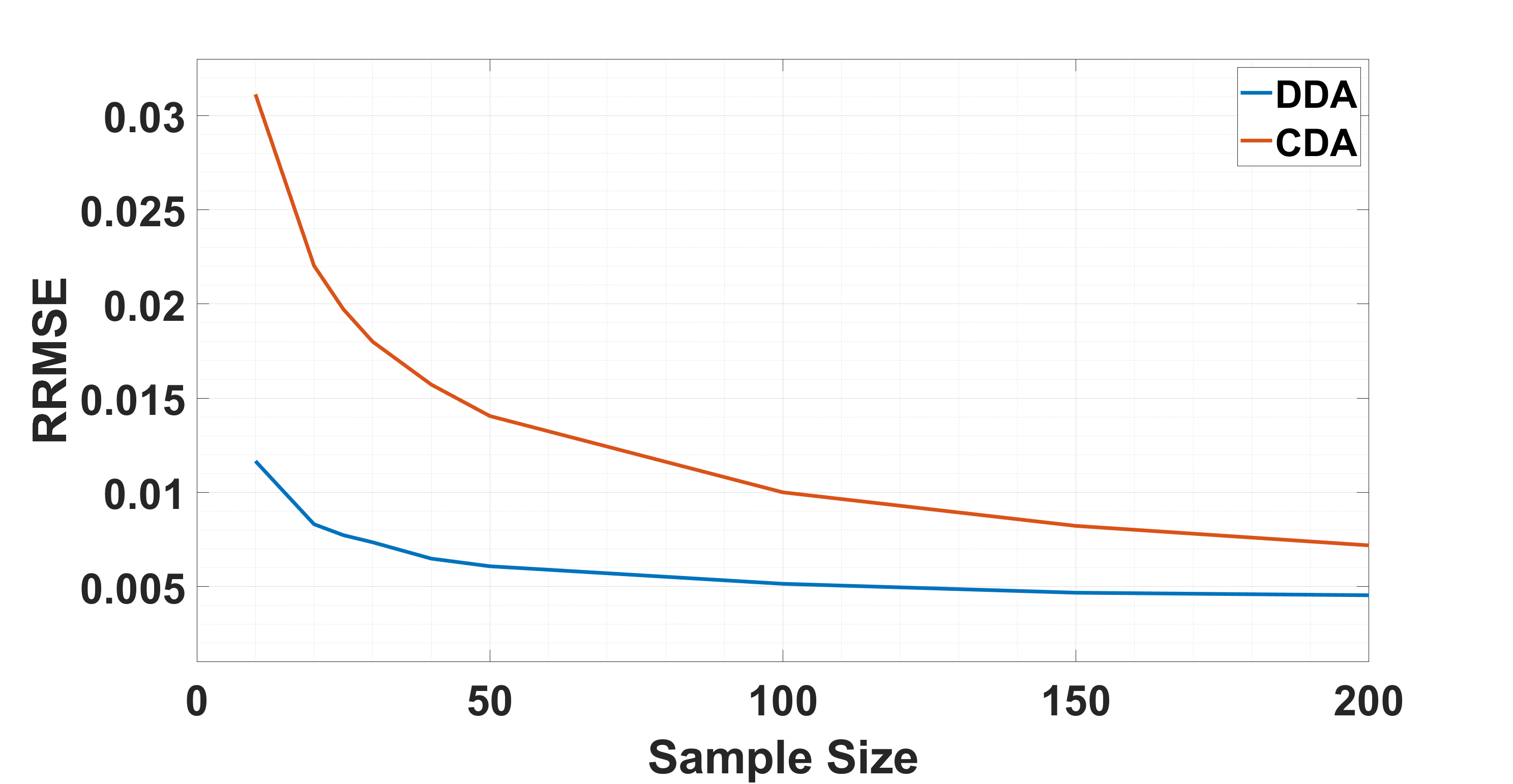}} \,


\caption{(a) Boxplots describing the distribution of the RRMSE of the downscaled temperature solutions for the DDA solution at the last time step (blue) and CDA solution at $t = 15$ (orange) for different sample sizes. 
(b) RRMSE curves of the average solution, at the final time, of the DDA (blue) and the CDA (orange) downscaled fields for different numbers of downscaled samples.}
\label{fig:expSol}
\end{figure}

Because the observational noise is unbiased, one may anticipate the mean solution to provide a better estimate of the reference solution than the individual realizations. 
The skill scores associated with the ensemble mean were computed for different subset sizes and compared to the distribution of the skill scores obtained from the $200$ nominal case solutions. 
Figure \ref{fig:expSol} illustrates the box plots of the asymptotic temperature RRMSE for both CDA and DDA downscaled fields, along with the curve of the RRMSE of the average solution for different numbers of realizations.
The results show that, for both the CDA and DDA and any number of realizations, the RRMSE of the average solution is lower than the RRMSE of any downscaled realization. 
The RRMSE of the average solution is also seen to decrease with larger number of realizations. 
Specifically, the RRMSE of the DDA downscaled solution dropped from approximately $3.5 \%$ to below $1 \%$ for more than $20$ realizations. 
Moreover, comparing the distributions from CDA to those from DDA, the results suggest that the DDA yields errors three times smaller than those of the CDA downscaled solutions. The CDA, however, shows a tighter distribution with a smaller spread between the dowscaled fields.

\subsection{Observing a Noisy Solution}
\label{sec:noisyObs}

In a practical setting, only a single noisy observation of the state is usually available.
Here, we generate samples of noisy observations starting by downscaling the given noisy coarse-scale observations to obtain a high-resolution downscaled reference.
Samples of noisy observations were then generated by adding iid noise sampled from a Gaussian distribution to this downscaled reference solution.
The resulting noisy observations are then considered as perturbed noisy observations. 
The sampled noise, at the nominal level, is added to the downscaled reference solution, as in Section \ref{res:noiseLevels}. 
By doing so, we investigate how CDA and DDA behave with the sampled noisy observations compared to the case when realizations of noisy observations of the truth are readily available.
Both CDA and DDA algorithms were then used to downscale these perturbed noisy observations for $\mathcal{S}=10$ and $\mathcal{R} = 15$.
The downscaled fields are compared to the reference fields to assess the performance of the downscaling algorithms.

\begin{figure}[!htbp]
\centering
\hspace{-1cm}

\subfloat{\includegraphics[width=130mm]{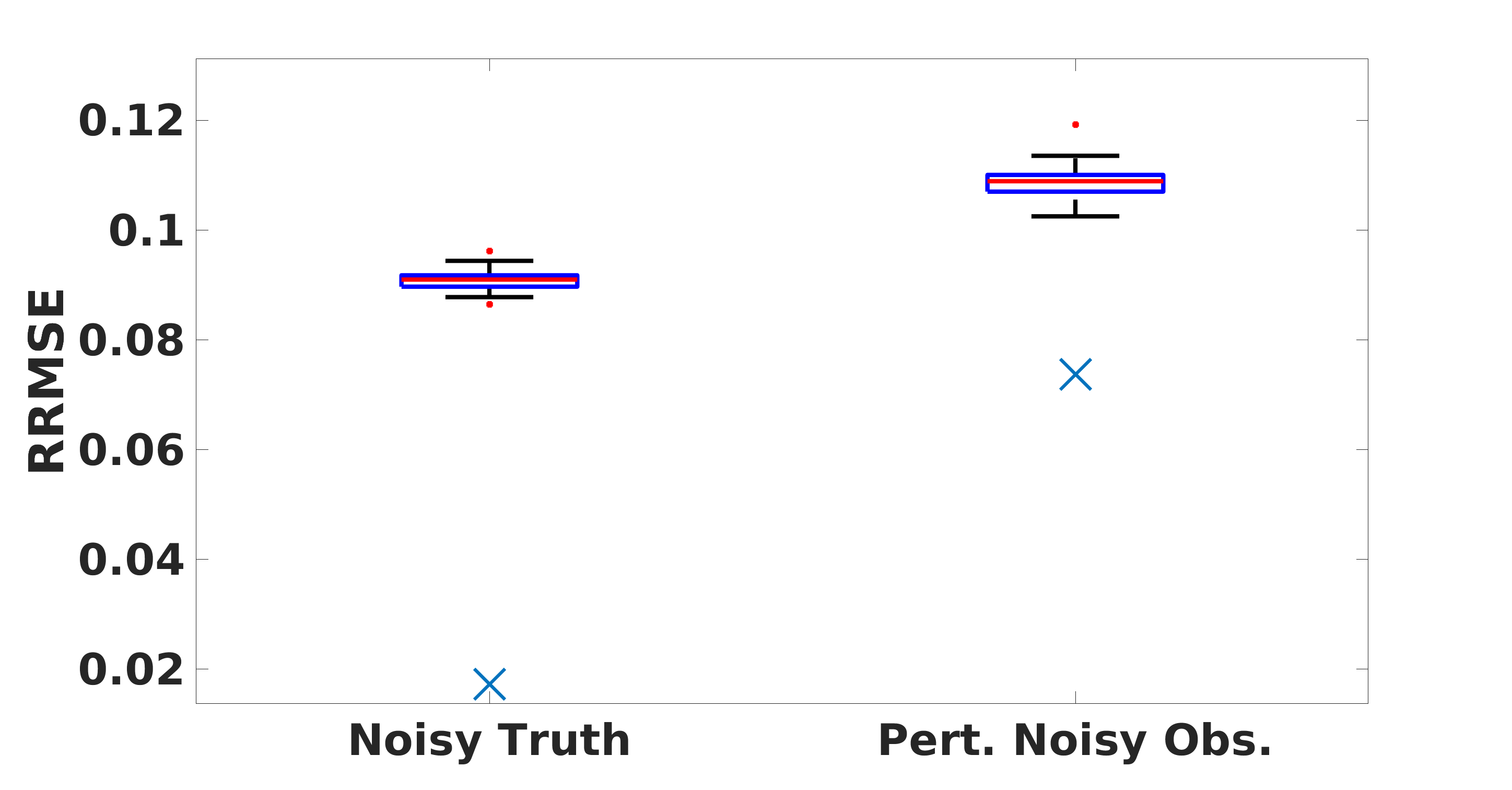}} \,

\caption{Box plots comparing the distributions of the RRMSE of the temperature sample solutions for the cases where observational data is provided based on the (left) truth and (right) noisy downscaled fields. The blue crosses represent the RRMSE of the average of all the downscaled samples.}
\label{fig:DoubleNoise}

\end{figure}

Figure \ref{fig:DoubleNoise} presents the box plots describing the distribution of the RRMSE of the DDA downscaled temperature fields relative to the reference temperature field at the final step. 
The box plots illustrate that the RRMSE increases by $20\%$ when the perturbed noisy observations were assimilated. 
Likewise, the variance of the downscaled solutions also increases and the RRMSE of the average solution is amplified by approximately $3.5$ times when downscaling was performed on the perturbed noisy observation. 
Nevertheless, the asymptotic RRMSE of the average solution was smaller than any of the asymptotic RRMSEs of the downscaled solutions, suggesting that it is always worth downscaling the perturbed noisy observations if computational resources are available.

\subsection{Relevance of Temperature Observations}
\label{sec:TempObsYAY_NAY}

The CDA was shown to converge to noise-free observational data using coarse-scale velocity observations alone \citep{FARHAT2015}.
The convergence of the algorithm, however, is faster when both temperature and velocity observations are assimilated \citep{Altaf2017}. 
Here, we assess the threshold for the temperature noise level beyond which temperature observations become disruptive to the convergence of CDA and DDA, and thus should be discarded.

In this experiment, observations of the velocity and temperature were provided for $\mathcal{S} = 15$ and $\mathcal{R} = 11$ with a fixed velocity noise level of $\sigma_{\mathbf{u}} = 0.05$ and three temperature noise levels: $\sigma_T = 0.05$, $0.1$ and $0.15$. 
An experiment was also conducted by downscaling using the coarse velocity observations only, which is equivalent to setting $\mu_{\Theta}=0$.

\begin{figure}[!htbp]
\centering
\hspace{-1cm}

\includegraphics[width = 0.9 \linewidth]{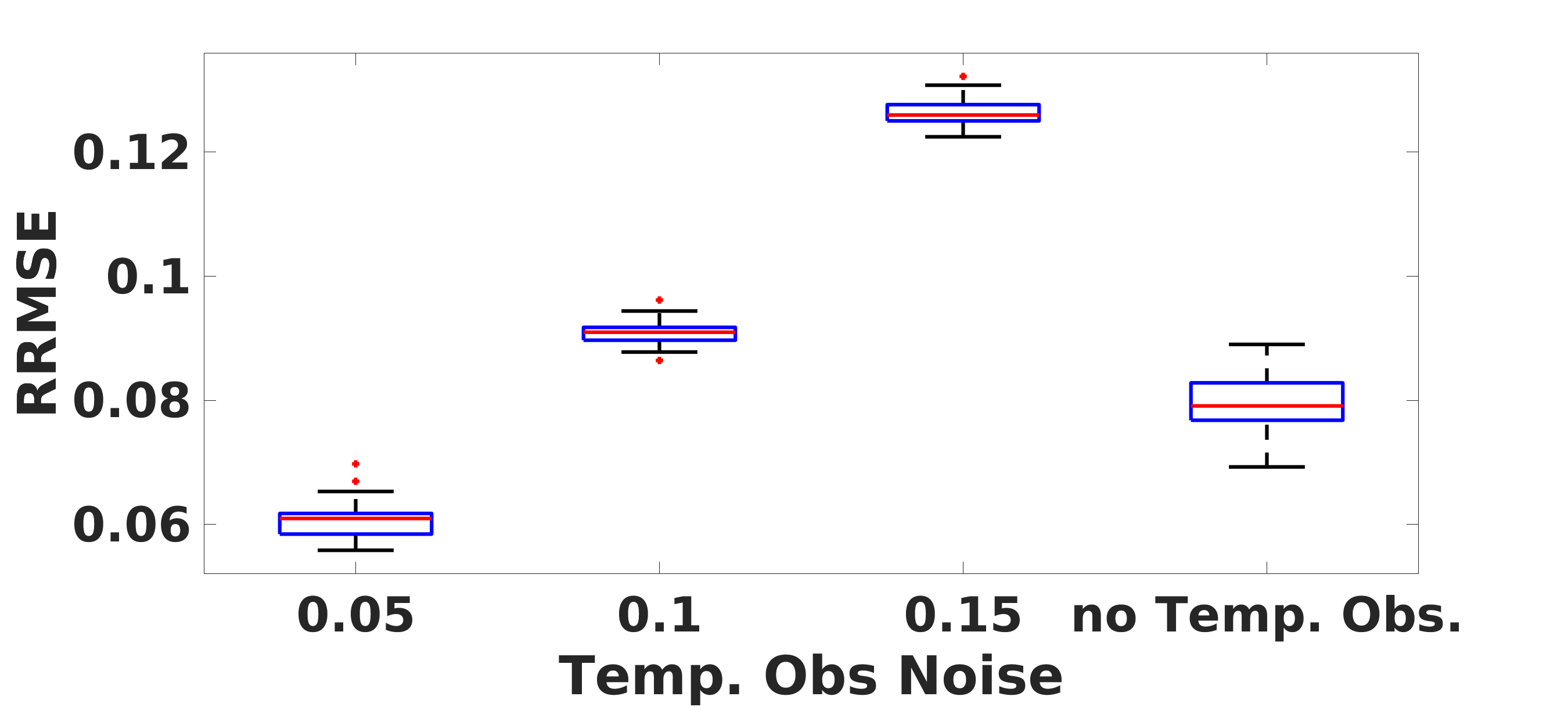}

\caption{Boxplot of the RRMSE of the DDA downscaled temperature solutions for temperature observations with low, medium, high temperature noise level, and no temperature observation. Each box plot is generated from a $50$ downscaled solutions at the last time step.}
\label{fig:lastStepTempObs}

\end{figure}

The time evolution of the RRMSE of the downscaled temperature fields illustrates that, in all tested cases, the RRMSE decreases until it reaches a plateau at approximately $t=40$. 
All RRMSE curves suggest similar characteristics before reaching the plateau, however, the spread is larger in the case when temperature was not assimilated. 
Figure \ref{fig:lastStepTempObs} presents the box plots describing the distribution of the asymptotic RRMSE of the DDA downscaled temperature. 
The downscaled solution without assimilating temperature observations achieves lower RRMSE than the case of high temperature noise level, comparable RRMSE with the case of medium temperature noise, and a larger RRMSE than that of the small temperature noise level. 
The case with no temperature observations shows the largest variance compared to all tested cases. 
This suggests that, in the case of medium to high noise, it might be better to discard temperature observations and opt for solely assimilating medium-noise velocity observations. 
On the other hand, if reliable temperature observations are available, it would be useful to jointly downscale with those of the velocity to obtain a more accurate solution.

\section{Conclusion}
\label{sec:conclusion}

This study assessed the performance of CDA, and a discrete variant called DDA, at downscaling the $2$D Rayleigh-B\'enard convection at various spatial and temporal resolutions of observations.
The analysis was conducted in an uncertain setting, where noisy observations of a reference solution are downscaled and statistically analyzed.
The CDA and DDA downscaled fields were contrasted qualitatively and quantitatively.
In particular, the performances of these algorithms were quantified by systematically analyzing the impact of observation noise level and the temporal and spatial resolutions of the observational grid, namely by quantifying the skill scores of the downscaled temperature and velocity fields and characterizing their statistics.

The skill scores of the downscaled fields from both CDA and DDA algorithms decay to a plateau, which reflects the convergence of the downscaled solution towards the reference field. 
Specifically, the CDA skill scores generally decay rapidly and monotonically, as opposed to the DDA skill scores, which exhibit a non-monotonic behavior before reaching a plateau. 
The plateaus of the DDA skill scores were at lower values compared to the CDA skill scores, meaning that errors of the DDA downscaled fields become smaller than those of CDA fields after the algorithm converges.

A thorough analysis of the functional dependence of the expected $\ell_2$-error of the downscaled fields ($\Lambda$) was performed on the noise level, temporal and spatial resolutions of the observation grid. 
$\Lambda$ was observed to vary quadratically with the variance of the observational noise for both CDA and DDA, in agreement with the theoretical results of \cite{Bessaih2015}. 
Moreover, the CDA results suggest that $\Lambda$ is proportional to the observation's temporal resolution, as opposed to a nonlinear dependence in the case of DDA. 
Specifically, for small and large temporal resolutions, $\Lambda$ is smaller when downscaling with CDA, however, for intermediate temporal resolutions, DDA downscaled fields result in a smaller $\Lambda$. 
Finally, our results illustrate that $\Lambda$ varies quadratically with the temporal resolution of the observational data in the case of CDA as opposed to the more sensitive DDA error where $\Lambda$ was observed to vary cubically with the temporal resolution of the observations.
The distribution of the downscaled fields were also examined by applying a statistical hypothesis test. 
Specifically, after reaching the plateau, the downscaled fields pass the Kolmogorov-Smirnov test at almost all tested locations in the domain, suggesting that they are normally distributed.

The precision of the estimated skill scores was also assessed following standard bootstrapping.
The results suggest that bootstrapping as few as $30$ realizations yields a distribution of the average RRMSE with a standard deviation comparable to that obtained from a $200$ realizations. 
Moreover, the mean downscaled solution achieves a lower error (roughly three times) than any of the downscaled solutions for both CDA and DDA.
When the temperature noise level is large, temperature observations degrade the quality of the downscaled fields and discarding them yields a smaller RRMSE. 
On the other hand, assimilating more accurate temperature observations helps reach the plateau faster and results in more truthful downscaled fields with smaller RRMSE.

Future work will extend the methodology to three-dimensional flows, address the combined effects of model and observation errors, and assess the performance of 
both CDA and DDA in comparison to other ensemble-based data assimilation techniques, such as the Ensemble Kalman Filter.
We further plan to perform a fundamental analysis of the behavior of the DDA algorithm, namely, to try to establish rigorous error estimates in settings involving observation and model noise.

\section*{Acknowledgments}

Research reported in this publication was supported by the Office of Sponsored Research (OSR) at King Abdullah University of Science and Technology (KAUST) CRG Award No. OSR-CRG2020-4336 and the Virtual Red Sea Initiative (Award No. REP/1/3268-01-01).

\section*{Data Availability}

The datasets generated during the current study are available from the corresponding author on reasonable request.


\bibliography{sn-bibliography}


\clearpage

\clearpage

\clearpage


\clearpage

\clearpage


\clearpage
\newpage

\clearpage
\newpage

\clearpage
\newpage


\clearpage
\newpage


\clearpage
\newpage






\clearpage
\newpage


\clearpage
\newpage


\clearpage
\newpage

\end{document}